\documentclass[12pt,twoside]{amsart}
\usepackage{amssymb}

\nonstopmode

\numberwithin{equation}{section}
\font\tengothic=eufm10 scaled\magstep 1
\font\sevengothic=eufm7 scaled\magstep 1
\newfam\gothicfam
          \textfont\gothicfam=\tengothic
          \scriptfont\gothicfam=\sevengothic

\newcommand{\Z}{\mathbb{Z}}

\DeclareMathOperator{\chara}{char}

\DeclareMathOperator{\codim}{codim}

\DeclareMathOperator{\im}{im}

\newcommand{\al}{\alpha}
\newcommand{\be}{\beta}
\newcommand{\ga}{\gamma}

\DeclareMathOperator{\pnt}{\raise 0.5mm \hbox{\large\bf.}}

\def\cocoa{{\hbox{\rm C\kern-.13em o\kern-.07em C\kern-.13em o\kern-.15em A}}}

\newcommand{\s}{\; | \;}

\newtheorem{theorem}{Theorem}[section]
\newtheorem{lemma}[theorem]{Lemma}
\newtheorem{proposition}[theorem]{Proposition}
\newtheorem{corollary}[theorem]{Corollary}
\newtheorem{conjecture}[theorem]{Conjecture}

\theoremstyle{definition}
\newtheorem{remark}[theorem]{Remark}
\newtheorem{example}[theorem]{Example}

\newtheorem{question}[theorem]{Question}
\newtheorem{problem}[theorem]{Problem}

\begin{document}

\title[Monomial ideals and the Weak Lefschetz Property]{Monomial
ideals, almost complete intersections and the Weak Lefschetz Property}


\author[J.\ Migliore, R.\ Mir\'o-Roig, U.\  Nagel]{Juan C.\ Migliore$^{*}$, Rosa M.\ Mir\'o-Roig$^{**}$, Uwe Nagel$^+$}
\address{Department of Mathematics, University of Notre Dame,
Notre Dame, IN 46556, USA}
\email{Juan.C.Migliore.1@nd.edu}

\address{Facultat de Matem\`atiques,
   Department d'Algebra i Geometria,
   Gran Via des les Corts Catalanes 585,
   08007 Barcelona,
   SPAIN }
\email{miro@ub.edu}

\address
{Department of Mathematics, University of Kentucky,
715 Patterson Office Tower, Lexington, KY 40506-0027, USA}
\email{uwenagel@ms.uky.edu}

\begin{abstract}
Many algebras are expected to have the Weak Lefschetz property
though this is often very difficult to establish. We illustrate the
subtlety of the problem by studying monomial and some closely
related ideals. Our results exemplify the intriguing dependence of
the property  on the characteristic of the ground field, and on
arithmetic properties of the exponent vectors of the monomials.
\end{abstract}

\thanks{${}^*$ Part of the work for this paper was done while the first
author was sponsored by the National Security Agency under Grant
Number H98230-07-1-0036.\\
${}^{**}$ Part of the work for this paper was done while the second
author was partially supported by MTM2007-61104.\\
${}^+$ Part of the work for this paper was done while the third
author was sponsored by the National Security Agency under Grant
Number H98230-07-1-0065.\\
The authors thank Fabrizio Zanello for useful and enjoyable
conversations related to some of this material.  They also thank
David Cook II for useful comments. \\
}

\maketitle

\section{Introduction}

Let $A$ be a standard graded Artinian algebra over the field $K$.
Then $A$ is said to have the {\em Weak Lefschetz property (WLP)}
if there is a linear form $L \in (A)_1$ such that, for all
integers $j$, the multiplication map
\[
\times L: (A)_{j-1} \to (A)_j
\]
has maximal rank, i.e.\ it is injective or surjective. In this
case, the linear form $L$ is called a {\em Lefschetz element} of
$A$.  (We will often abuse notation and say that the corresponding
ideal has the WLP.)  The Lefschetz elements of $A$ form a Zariski
open, possibly empty, subset of $(A)_1$. Part of the great
interest in the WLP stems from the fact that its presence puts
severe constraints on the possible Hilbert functions (see
\cite{HMNW}), which can appear in various disguises (see, e.g.,
\cite{St-faces}). Though many algebras are expected to have the
WLP, establishing this property is often rather difficult. For
example, it is open whether every complete intersection of height
four over a field of characteristic zero has the WLP. (This is
true if the height is at most 3 by \cite{HMNW}.)

In some sense, this note presents a case study of the WLP for
monomial ideals and almost complete intersections. Our results
illustrate how subtle the WLP is. In particular, we investigate
its dependence on the characteristic of the ground field $K$.  The
following example (Example \ref{both cases}) illustrates the
surprising effect that the characteristic can have on the WLP.
Consider the ideal $I = (x^{10}, y^{10},z^{10}, x^3 y^3 z^3)
\subset R = K[x,y,z]$.  Our methods show that $R/I$ fails to have
the WLP in characteristics 2, 3 and 11, but possesses it in all
other characteristics.

One starting point of this paper has been Example 3.1 in
\cite{BK}, where Brenner and Kaid show that, over an algebraically
closed field of characteristic zero, any  ideal of the form
$(x^3,y^3,z^3,f(x,y,z))$, with $\deg f = 3$, fails to have the WLP
if and only if $f \in (x^3,y^3,z^3,xyz)$.  In particular, the
latter ideal is the {\em only} such monomial ideal that fails to
have the WLP. This paper continues the study of this question.

The example of Brenner and Kaid satisfies several  interesting
properties.  In this paper we isolate several of these properties
and examine the question of whether or not the WLP holds for such
algebras, and we see to what extent we can generalize these
properties and still get meaningful results.   Some of our results
hold over a field of arbitrary characteristic, while others show
different ways in which the characteristic plays a central role in
the WLP question.  (Almost none are characteristic zero results.)
Most of our results concern monomial ideals, although in Section
\ref{almost monomial} and Section \ref{final comments} we show
that even minor deviations from this property can have drastic
effects on the WLP.  Most of our results deal with almost complete
intersections in three or more variables, but we also study ideals
with more generators (generalizing that of Brenner and Kaid in a
different way).

More specifically, we begin in  Section \ref{sec:tools-for-WLP} with
some simplifying tools for studying the WLP.   These are applied
throughout the paper.  We also recall the construction of basic
double linkage.

In Section \ref{class of mon ideals} we consider the class of
monomial ideals in $K[x_1,\ldots,x_r]$ of the form
\[
(x_1^k, x_2^k,\dots,x_r^k) + (\hbox{all squarefree monomials of degree $d$}).
\]
Note that the example of Brenner and Kaid is of this form.  Our
main result in this section (Theorem \ref{thm-wlp}) says that when
$d=2$ we always have the WLP, but if $d=3$ and $k \geq 2$ then we
have two cases: if $K$ has characteristic 2 then we never have the
WLP, but if the characteristic is not 2 then we have the WLP if
and only if $k$ is even.

In Section \ref{arb codim},  we consider almost complete
intersections of the form $(x_1^r,\dots,x_r^r, x_1\cdots x_r )$
with $r \geq 3$ (note that the result of Brenner and Kaid dealt
with the case $r=3$ in characteristic zero).  Our main result for
these algebras is that they always fail to have the WLP,
regardless of the characteristic.  The proof is surprisingly
difficult.

In Section \ref{almost monomial} we  explicitly illustrate the
fact that  even a minuscule change in the ideal can affect the
WLP. Specifically, we consider the ideals of the form
\[
(x_1^r,\dots,x_r^r, x_1\cdots x_{r-1} \cdot (x_1+x_r) ).
\]
 We show that this has
the same Hilbert function as the corresponding ideal in the previous
section, but the WLP behavior is very different. For example, the
two ideals
\[
(x_1^4,\ldots,x_4^4, x_1 x_2 x_3 x_4) \quad \text{and} \quad
(x_1^4,\ldots,x_4^4, x_1 x_2 x_3 (x_1 + x_4))
\]
have the same Hilbert function, but the former never has the WLP
while the latter has the WLP if and only if the characteristic of
$K$ is not two or five.

In Section \ref{monomial aci} we  turn to monomial almost complete
intersections in three variables, generalizing the Brenner-Kaid
example in a different direction.  To facilitate this study, we
assume that the algebra is also level (as is the case for Brenner
and Kaid's example).  We give a number of results  in this section,
which depend on the exponent vectors of the monomials. We end with a
conjectured classification of the level Artinian monomial ideals in
three variables that fail to have the WLP (Conjecture \ref{level wlp
conj}). The work in  Sections \ref{monomial aci} and \ref{half conj}
proves most of this conjecture. We end the paper in Section
\ref{final comments} with some suggestive computations and natural
questions coming from our work.

\tableofcontents


\section{Tools for studying the WLP}
\label{sec:tools-for-WLP}

In this section we establish various general results that help to
study the WLP and that are used throughout the remainder of this
paper.  Throughout this paper we set $R = K[x_1,\dots,x_r]$, where $K$ is a field.  Sometimes we will have specific values of $r$ (usually 3) and sometimes we will have further restrictions on the field $K$.

Our first results singles out the crucial maps to be studied if we
consider the WLP of a level algebra. Recall that an
Artinian algebra
is called {\em level} if its socle is concentrated in one degree.

\begin{proposition} \label{gen wlp}
Let $R/I$ be an Artinian standard graded algebra and let $L$ be a general linear form.  Consider the homomorphisms $\phi_d : (R/I)_d \rightarrow (R/I)_{d+1}$ defined by multiplication by $L$, for $d \geq 0$.

\begin{itemize}
\item[(a)] If $\phi_{d_0}$ is surjective for some $d_0$ then $\phi_d$ is surjective for all $d \geq d_0$.

\item[(b)] If $R/I$ is level and $\phi_{d_0}$ is injective for some $d_0 \geq 0$ then $\phi_d$ is injective for all $d \leq d_0$.

\item[(c)] In particular, if $R/I$ is level and $\dim (R/I)_{d_0}
= \dim (R/I)_{d_0+1}$ for some $d_0$ then $R/I$ has the WLP if and
only if $\phi_{d_0}$ is injective (and hence is an isomorphism).
\end{itemize}
\end{proposition}

\begin{proof}
Consider the exact sequence
\[
0 \rightarrow \frac{[I:L]}{I} \rightarrow R/I \stackrel{\times L}{\longrightarrow} (R/I) (1) \rightarrow (R/(I,L))(1) \rightarrow 0
\]
where $\times L$ in degree $d$ is just $\phi_d$.  This shows that the cokernel of $\phi_d$ is just $(R/(I,L))_{d+1}$ for any $d$. If $\phi_{d_0}$ is surjective, then $(R/(I,L))_{d_0 +1} = 0$, and the same necessarily holds for all subsequent twists since $R/I$ is a standard graded algebra.  Then (a) follows immediately.

For (b), recall that the $K$-dual of the finite length module $R/I$ is a shift of the canonical module of $R/I$, which we will denote simply by $M$.  Since $R/I$ is level, $M$ is generated in the first degree.  But now if we consider the graded homomorphism of $M$ to itself induced by multiplication by $L$, a similar analysis (recalling that $M$ is generated in the first degree) gives that once this multiplication is surjective in some degree, it is surjective thereafter.  The result on $R/I$ follows by duality.

Part (c) follows immediately from (a) and (b).
\end{proof}

If the field is infinite and the $K$-algebra satisfies the WLP for some linear form, then it does for a general linear form. For monomial ideals there is no need to consider a general linear form.

\begin{proposition}
  \label{lem-L-element}
Let $I \subset R$ be an Artinian monomial ideal and assume that the field $K$ is infinite. Then $R/I$ has the WLP if and only if  $x_1 + \cdots + x_r$ is a Lefschetz element for $R/I$.
\end{proposition}

\begin{proof}
Set $A = R/I$ and let $L = a_1 x_1 + \cdots + a_r x_r$ be a general linear form in $R$. Thus, we may assume that each coefficient $a_i$ is not zero and, in particular, $a_r = 1$. Let $J \subset S := K[x_1,\ldots,x_{r-1}]$ be the ideal that is generated by elements that are obtained from the minimal generators of $I$ after substituting $a_1 x_1 + \cdots + a_{r-1} x_{r-1}$ for $x_r$. Then $A/L A \cong S/J$.

Each minimal generator of $J$ is of the form $x_1^{j_1} \cdots
x_{r-1}^{j_{r-1}} (a_1 x_1 + \cdots + a_{r-1} x_{r-1})^{j_r}$.
Replacing it by $(a_1 x_1)^{j_1} \cdots (a_{r-1} x_{r-1})^{j_{r-1}}
(-a_1 x_1 + \cdots - a_{r-1} x_{r-1})^{j_r}$ does not change the
ideal $J$ because $a_1 \cdots a_{r-1} \neq 0$. Using the isomorphism
$K[y_1,\ldots,y_{r-1}] \to S, \ y_i \mapsto a_i x_i$, we see that
$A/L A$ and  $A/(x_1 + \cdots + x_r) A$ have the same Hilbert
function. Since we can decide whether $L$ is a Lefschetz element for
$A$ by solely looking at the Hilbert function of $A/L A$, the claim
follows.
\end{proof}

If $A$ is an Artinian $K$-algebra with the WLP and $E$ is an
extension field of $K$, then also $A \otimes_K E$ has the WLP.
However, the converse is not clear. We pose this as a problem.

\begin{problem}
  \label{prop-ext}
Is it true that  $A$ has the WLP if and only if $A \otimes_K E$ has the WLP?
\end{problem}

Proposition \ref{lem-L-element} shows that the answer is affirmative
in the case of monomial ideals.

\begin{corollary}
  \label{cor-WLP-inv}
Let $E$ be an extension field of the infinite field $K$.  If $I
\subset R$ is an Artinian monomial ideal, then $R/I$ has the WLP if
and only if $(R/I) \otimes_K E$ does.
\end{corollary}

The following result applies if we can hope that the multiplication
by a linear form is surjective.

\begin{proposition}
\label{splitting and surj}
Let $I \subset R = K[x_1,\dots,x_r]$, where $K$
is a field and $A = R/I$ is Artinian.  Let $d$ be any
degree such that $h_A(d-1) \geq h_A(d)
>0$. Let $L$ be a linear form,  let $\bar R = R/(L)$ and let
$\bar I$ be the image of $I$ in $\bar R$. Denote by $\bar A$ the
quotient $\bar R / \bar I$. Consider the minimal free $\bar
R$-resolution of $\bar A$:
\[
0 \rightarrow \bigoplus_{i=1}^{p_{r-1}} \bar R(-b_i) \rightarrow \cdots \rightarrow \bigoplus_{j=1}^{p_1} \bar R(-a_j) \rightarrow \bar R \rightarrow  \bar A \rightarrow 0
\]
where $a_1 \leq \dots \leq a_{p_1}$ and  $b_1 \leq \dots \leq
b_{p_{r-1}}$.  Then the following are equivalent:

\begin{itemize}

\item[(a)] the multiplication by $L$ from $A_{d-1}$ to
$A_d$ fails to be surjective;

\item[(b)] $\bar A_d \neq 0$;

\item[(c)]  $b_{p_{r-1}} \geq d+r-1$;

\item[(d)]
Let $G_1,\dots,G_{r-1}$ be a  regular sequence in $\bar I$ of
degrees $c_1,\dots,c_{r-1}$ respectively that extends to a
minimal generating set for $\bar I$.  Then there exists a form
$F \in \bar R$ of degree $\leq c_1+ \cdots +c_{r-1} - (d+r-1)$,
non-zero modulo $(G_1,\dots,G_{r-1})$, such that $F \cdot \bar I
\subset (G_1,\dots,G_{r-1})$.
\end{itemize}
\end{proposition}

\begin{proof}
From the exact sequence
\[
\dots \rightarrow A_{d-1} \stackrel{\times L}{\longrightarrow} A_d \rightarrow (R/(I,L))_d \rightarrow 0
\]
it follows that the multiplication fails to be surjective if and only if $\bar A_d = (R/(I,L))_d \neq 0$.  The latter holds if and only if
\[
d \leq \hbox{socle degree of } \bar A = b_{p_{r-1}}  -(r-1),
\]
from which the equivalence of (a), (b) and (c) follows.

To show the equivalence of (c) and (d) we invoke liaison theory.  Let
\[
J = (G_1,\dots,G_{r-1}) : \bar I.
\]
  A free resolution for $J$ can be obtained from that of $\bar I$ and $(G,G')$ by a standard mapping cone argument (see for instance \cite{migbook}), as follows.  We have the following commutative diagram (where the second one is the Koszul resolution for $(G_1,\dots,G_{r-1})$):
\[
\begin{array}{cccccccccccccccccccccc}
0 & \rightarrow & \bigoplus_{i=1}^{p_{r-1}} \bar R(-b_i) & \rightarrow & \cdots & \rightarrow & \bigoplus_{j=1}^{p_1} \bar R(-a_j)  & \rightarrow & \bar I &  \rightarrow & 0 \\
&& \uparrow &&&& \uparrow && \uparrow \\
0 & \rightarrow & \bar R(-c_1-\cdots -c_{r-1}) & \rightarrow & \cdots & \rightarrow & \bigoplus_{k = 1}^{r-1} \bar R(-c_k)  & \rightarrow & (G,G') & \rightarrow & 0
\end{array}
\]
where the rightmost vertical arrow is an inclusion.  This yields a free resolution for $J$ (after splitting $ \bigoplus_{k=1}^{r-1} \bar R(-c_k)$ and re-numbering the $a_j$, and setting $c := c_1 + \cdots + c_{r-1}$):
\[
0 \rightarrow \bigoplus_{j=1}^{p_1 -(r-1)} \bar R(a_j - c) \rightarrow \cdots \rightarrow \bigoplus_{i=1}^{p_{r-1}} \bar R(b_i - c) \oplus \bigoplus_{k=1}^{r-1} \bar R(-c_k) \rightarrow J \rightarrow 0.
\]
Clearly $b_{p_{r-1}} \geq d+r-1$ if and only if $J$ has a minimal generator, $F$, of degree $\leq c-(d+r-1)$.    The result then follows from the definition of $J$ as an ideal quotient.
\end{proof}

We conclude this section by recalling a concept from liaison theory,
which we do not state in the greatest generality.

Let $J \subset I \subset R = K[x_1,\ldots,x_r]$ be homogeneous
ideals such that $\codim J = \codim I - 1$. Let $\ell \in R$ be a
linear form such that $J  : \ell = J$. Then the ideal $I' := \ell
\cdot I + J$ is called a {\em basic double link} of $I$. The name
stems from the fact that $I'$ can be Gorenstein linked to $I$ in two
steps if $I$ is unmixed and $R/J$ is Cohen-Macaulay and generically
Gorenstein (\cite{KMMNP}, Proposition 5.10). However, here we only
need the relation among the Hilbert functions.

\begin{lemma}
  \label{lem:BDL}
For each integer $j$,
\[
\dim_K (R/I')_j = \dim_K (R/I)_{j-1} + \dim_K (R/J)_{j} - \dim_K
(R/J)_{j-1}.
\]
\end{lemma}

\begin{proof}
  This follows from the exact sequence (see \cite{KMMNP}, Lemma 4.8)
\[
0 \to J(-1) \to J  \oplus I (-1) \to I' \to 0.
\]
\end{proof}


\section{A class of monomial ideals}
\label{class of mon ideals}

We now begin our study of a certain class of Artinian monomial
ideals.
Let $I_{r,k,d}$ be the monomial ideal defined by
\begin{equation}
  \label{eq:def-ideal}
(x_1^k, x_2^k,\dots,x_r^k) + (\hbox{all squarefree monomials of degree $d$}).
\end{equation}


Our first observation follows immediately be determining the socle
of $R/I_{r,k,d}$. It shows that we may apply Proposition \ref{gen
wlp}.

\begin{proposition} \label{inv sys}
The inverse system for $I_{r,k,d}$ is generated by the module generated by all monomials of the form
\[
x_{i_1}^{k-1} \cdots x_{i_{d-1}}^{k-1}.
\]
\end{proposition}

\begin{corollary}
The algebra $R/I_{r,k,d}$ is level of socle degree $(k-1)(d-1)$ and
socle type $\binom{r}{d-1}$.
\end{corollary}

Concerning the WLP we have:

\begin{theorem}
  \label{thm-wlp}
Consider the ring $R/I_{r,k,d}$.
\begin{itemize}
  \item[(a)] If $d=2$, then it has the WLP.
\item[(b)] Let $d = 3$ and $k \geq 2$.  Then:
\begin{itemize}
  \item[(i)] If $K$ has characteristic two, then $R/I_{r,k,d}$ does not have the WLP.
\item[(ii)] If the characteristic of $K$ is not two, then $R/I_{r,k,d}$ has the WLP
if and only if $k$ is even.
\end{itemize}
\end{itemize}
\end{theorem}

\begin{proof}
  For simplicity, write $I = I_{r,k,d}$ and $A = R/I_{r,k,d}$.

Claim (a) follows easily from the observation that $A$ has socle degree $k-1$ and that up to degree $k-1$ the ideal $I$ is radical, so multiplication by a general linear form is injective in degree $\leq k-1$.

To show claim (b) we first describe bases   of $(A)_{k-1}$ and
$(A)_k$, respectively. We choose the residue classes of the
elements in the following two sets.
\begin{eqnarray*}
  B_{k-1} & = & \{x_i^j x_m^{k-1-j} \s 1 \leq i < m \leq r,\; 1 \leq j \leq k-2 \} \cup
\{x_i^{k-1} \s 1 \leq i \leq r\} \\
B_k & = & \{x_i^j x_m^{k-j} \s 1 \leq i < m \leq r,\; 1 \leq j  \leq k-1 \}.
\end{eqnarray*}
Counting we get
\begin{equation}
\label{eq-hilb}
  h_A (k-1) = (k - 2) \binom{r}{2} + r \leq (k-1) \binom{r}{2},
\end{equation}
where the inequality follows from $r \geq d = 3$.

Now we assume that $k$ is odd. In this case we claim that $A$ does not have the WLP.
 Because of Inequality (\ref{eq-hilb}), this follows once we have shown that, for each linear form $L \in R$, the multiplication map
 $\phi_k: (A)_{k-1}  \stackrel{\times L}{ \longrightarrow} (A)_k$ is not injective.

To show the latter assertion we exhibit a non-trivial element in its kernel.
Write $L = a_1 x_1 + \ldots + a_r x_r$ for some $a_1,\ldots,a_r \in K$.
We define the polynomial $f \in R$ as
\[
f = \sum_{ j_i + j_m = k-1 } (-1)^{\max \{j_i, j_m\}}
(a_i x_i)^{j_i} (a_m x_m)^{j_m}.
\]
Note that $f$ is not in $I$.
We claim that $L \cdot f$ is in $I$. Indeed, since all monomials involving three distinct variables are in $I$, a typical monomial in $L \cdot f \!\!\mod I$ is of the form
\[
(a_i x_i)^{j_i} (a_m x_m)^{k-j_i}.
\]
It arises in exactly two ways in $L f$, namely as $(a_i x_i) \cdot (a_i x_i)^{j_i - 1} (a_m x_m)^{k-j_i}$ and as $(a_m x_m) \cdot (a_i x_i)^{j_i} (a_m x_m)^{k-1-j_i}$. Using that $k-1$ is even, it is easy to see that these two monomials occur in $f$ with different signs. It follows that the above multiplication map is not injective.

If $k$ is even, but $\chara K = 2$, then the same analysis again shows that $\phi_k$ is not injective. Hence, for the remainder of the proof we may assume that the characteristic of $K$ is not two.

Assume $k$ is even. Then we claim that $L = x_1 + \cdots + x_r$ is
a Lefschetz element. To this end we first show that the
multiplication map $\phi_{k}: (A)_{k-1}  \stackrel{\times L}{
\longrightarrow} (A)_{k}$ is injective.

 Let $f$ be any element in the vector space generated by $B_{k-1}$. Pick three of the variables $x_1,\ldots,x_r$ and call them $x, y, z$. Below we explicitly list all the terms in $f$ that involve only the variables $x, y, z$:
\begin{eqnarray*}
f & = & a_0 x^{k-1} + a_1 x^{k-2} y + \cdots + a_{k-2} x y^{k-2} + \\
&& \hspace*{1.7cm} b_1 x^{k-2} z + \cdots + b_{k-2} x z^{k-2} + b_{k-1} z^{k-1} + \\
&&c_0 y^{k-1} + c_1 y^{k-2} z + \cdots + c_{k-2} y z^{k-2} \\
&& + \cdots.
\end{eqnarray*}
As above, we see that each monomial in $L \cdot f$ arises from exactly two of the monomials in $f$. Hence the condition $L \cdot f \in I$ leads to the following three systems of equations.  Focussing  only on the variables $x, y$ we get:
\begin{eqnarray*}
  \label{eq-sys1}
a_0 + a_1 & = & 0 \\
a_1 + a_2 & = & 0 \\
 & \vdots & \\
a_{k-3} + a_{k-2} & = &0 \\
a_{k-2} + c_0 & = & 0.
\end{eqnarray*}
It follows that $a_i = (-1)^i a_0$ and
\begin{equation}
  \label{eq-1}
c_0 = (-1)^{k-2} a_0 = a_0
\end{equation}
because $k$ is even. Considering the variables $x, z$ we obtain:
\begin{eqnarray*}
  \label{eq-sys2}
a_0 + b_1 & = & 0 \\
b_1 + b_2 & = & 0 \\
 & \vdots & \\
b_{k-2} + b_{k-1} & = & 0,
\end{eqnarray*}
hence
\begin{equation}
  \label{eq-2}
b_i = (-1)^i a_0.
\end{equation}
Finally, using the variables $y, z$ we get:
\begin{eqnarray*}
  \label{eq-sys3}
c_0 + c_1 & = & 0 \\
 & \vdots & \\
c_{k-1} + c_{k-2} & = & 0 \\
c_{k-2} + b_{k-1} & = & 0.
\end{eqnarray*}
Combining this, it follows that
\[
-a_0 = b_{k-1} = -c_0 = a_0.
\]
Since we assumed that the characteristic of $K$ is not two, we conclude that the three linear systems above have only the trivial solution. Since the variables $x, y, z$ were chosen arbitrarily, we see that the map $\phi_k$ is injective, as claimed.

According to Lemma \ref{gen wlp} it remains to show that the
multiplication map $\phi_{k+1}: (A)_{k}  \stackrel{\times L}{
\longrightarrow} (A)_{k+1}$ is surjective. Note that the residue
classes of the elements of the form $x_i^j x_m^{k+1-j}$ with $2
\leq j \leq k-1,\; 1 \leq i < m \leq r$ form a basis of
$(A)_{k+1}$. Setting for simplicity $x := x_i, y = x_m$ it is
enough to show that, for each $j = 2,\ldots,k-1$, the residue
class of $x^j y^{k+1-j}$ is in the image of $\phi_{k+1}$.

We induct on $j \geq 2$. If $j=2$, then we get modulo $I$ that $L \cdot x y^{k-1} \equiv x^2 y^{k-1}$, thus $\overline{x^2 y^{k-1}} \in \im \phi_{k+1}$, as claimed. Let $3 \leq j \leq k-1$, then, modulo $I$, we get $L \cdot x^{j-1} y^{k-j} \equiv x^j y^{k-j} + x^{j-1} y^{k-j+1}$. Since by induction $\overline{x^{j-1} y^{k-j+1}} \in \im \phi_{k+1}$, we also obtain $\overline{x^{j} y^{k-j}} \in \im \phi_{k+1}$. This completes the proof.
\end{proof}

The above result and our computer experiments suggest that the larger $d$ becomes, the rarer it is that $R/I_{r,k,d}$ has  the WLP.
 Based on computer experiments we expect the following to be true.

\begin{conjecture}
  \label{conj:wlp}
Consider the algebra $R/I_{r,k,d}$. Then
\begin{itemize}
\item[(a)] If $d = 4$, then it has the WLP if and only if $k\!\!\!\mod 4$ is 2 or 3.

\item[(b)] If $d = 5$, then the WLP fails.

\item[(c)] If $ d= 6$, then the WLP fails.
\end{itemize}

\end{conjecture}

We summarize our results in case $k = d= 3$.

\begin{example} \label{1 r rchoose2}
Consider the ideal
\[
I_{r,3,3} = (x_1^3,x_2^3,\dots,x_r^3, \hbox{\rm (all squarefree  monomials of degree 3)}).
\]
Then the corresponding inverse system is  $(x_1^2x_2^2, \ \
x_1^2x_3^2, \ \ \dots, \ \ x_{r-1}^2 x_r^2 )$.  Furthermore, the
Hilbert function of $R/I_{r,3,3}$ is
\[
1 \ \ r \ \ \binom{r+1}{2} \ \ r(r-1) \ \ \binom{r}{2} \ \ 0
\]
and $R/I_{r,3,3}$ fails to have the WLP because the map from
degree 2 to degree 3 by a general linear form is not injective.
\end{example}

\begin{remark}
By truncating, we get a compressed level algebra with socle degree
3 that fails to have the WLP.  We expect that there are compressed
level algebras with larger socle degree that fail to have the WLP.
However, we do not know such an example.
\end{remark}



\section{Monomial almost complete intersections in any codimension}
\label{arb codim}


In the paper \cite{MMR3} the first and second authors asked the following question (Question 4.2, page 95): For any integer $n \geq 3$, find the minimum number $A(n)$ (if it exists) such that {\em every} Artinian ideal
$I \subset K[x_1,\dots,x_n]$ with number of generators $\mu (I) \leq A(n)$ has the WLP.
  In Example \ref{A(n) example} below, we show that $A(n)$ does not exist in positive characteristic.  In any case, in \cite{BK} it was shown for $n=3$ and characteristic zero that $A(3) = 3$ (also using a result of \cite{HMNW}), as noted in the introduction.  A consequence of the main result of this section, below, is that in any number of variables and any characteristic there is an almost complete intersection that
   fails to have the
  WLP.  Hence the main open question that remains is whether, in characteristic zero,
  all complete intersections have the WLP (as was shown for $n=3$ in \cite{HMNW}), i.e.\ whether $A(n) = n$ in characteristic zero.

We begin by considering ideals of the form
\begin{equation} \label{defofIrk}
I_{r,k} = (x_1^k,\dots,x_r^k ,x_1 \dots x_r) \subset K[x_1,\dots,x_r].
\end{equation}
Notice that this is a special case of the class of ideals  described
in Section \ref{class of mon ideals}. It is not too difficult to
determine the graded Betti numbers.

\begin{proposition}
The minimal free resolution of $I_{r,k}$ has the form
{\footnotesize
\[
\begin{array}{c}
0 \rightarrow
\begin{array}{c}
 R(-r+(r-1)(1-k))^{\binom{r}{r-1}}
 \end{array}

 \rightarrow

 \begin{array}{c}
 R(-(r-1)k))^{\binom{r}{r-1}} \\
 \oplus \\
 R(-r + (r-2)(1-k))^{\binom{r}{r-2}}
 \end{array}

\rightarrow \cdots \rightarrow

\begin{array}{c}
R(-3k)^{\binom{r}{3}} \\
\oplus \\
R(-r+2(1-k))^{\binom{r}{2}}
\end{array}

\rightarrow

\\ \\ \\

\begin{array}{c}
R(-2k)^{\binom{r}{2}} \\
\oplus \\
R(-r+(1-k))^r
\end{array}

\rightarrow

\begin{array}{c}
R(-k)^r \\
\oplus \\
R(-r)
\end{array}
\rightarrow I_{r,k} \rightarrow 0.

\end{array}
\]
}

\end{proposition}

\begin{proof}
Since $I_{r,k}$ is an almost complete  intersection, we can link it
using the complete intersection ${\mathfrak a} =(x_1^k,\dots,x_r^k)$
to an Artinian Gorenstein ideal, $J$.  However, since both $I_{r,k}$
and $\mathfrak a$ are monomial, so is $J$.  But it was first shown
by Beintema \cite{beintema} that any monomial Artinian Gorenstein
ideal is a complete intersection.  Hence we get by direct
computation that $(x_1^k,\dots,x_r^k) : x_1x_2\cdots x_r =
(x_1^{k-1},\dots,x_r^{k-1})$.  Then use the mapping cone and observe
that there is no splitting.
\end{proof}

%
%
%

Before we come to the main result of this section,  we prove a
preliminary result about the Hilbert function of  complete
intersections that will allow us to apply Proposition \ref{splitting
and surj}.

\begin{lemma} \label{diff of hf}
Let $R = K[x_1,\dots,x_s]$ with $s \geq 2$, and let
\[
I_s = (x_1^s,\dots,x_s^s) \hbox{\ \ \ \ and \ \ \ \ } J_s = (x_1^{s+1},x_2^{s+1},x_3^s,\dots,x_s^s).
\]
Note that the midpoint of the Hilbert function of $R/I_s$ is $\binom{s}{2}$ and that of $R/J_s$ is $\binom{s}{2} +1$.  Then
\[
\hbox{$h_{R/I_s}\left ( \binom{s}{2} \right ) - h_{R/I_s} \left ( \binom{s}{2} -1 \right )
\leq
h_{R/J_s} \left ( \binom{s}{2} +1 \right ) - h_{R/J_s} \left ( \binom{s}{2} + 2 \right ).$}
\]
\end{lemma}

\begin{proof}
The lemma is trivial to verify when $s = 2$ or $s = 3$, so we assume
$s \geq 4$ for this proof. Observe that both quantities are
positive, but one  involves a difference to the left of the midpoint
of the Hilbert function, while the other involves a difference to
the right.  We will use this formulation, although there exists
others thanks to the symmetry of the Hilbert function of an Artinian
complete intersection.

Our approach will be via basic double linkage. We will use the
formula in Lemma \ref{lem:BDL}  without comment.   In fact, $J_s$ is
obtained from $I_s$ by a sequence of two basic double links:
\[
\begin{array}{ll}
I_s & \leadsto x_1 \cdot I_s + (x_2^s,\dots,x_s^s) := G = (x_1^{s+1},x_2^s,\dots,x_s^s) \\
& \leadsto x_2 \cdot G + (x_1^{s+1},x_3^s,\dots,x_s^s) = J_s.
\end{array}
\]
Note that $G$ is a complete intersection of codimension $s$ and
that the ideals $C_1 := (x_2^s,\dots,x_s^s)$ and $C_2 :=
(x_1^{s+1},x_3^s,\dots,x_s^s)$ are complete intersections of
codimension $s-1$.  The midpoints of the {\em $h$-vectors } of
$R/C_1$ and $R/C_2$ are $\frac{(s-1)^2}{2}$ and
$\frac{(s-1)^2+1}{2}$ respectively.  We now compute Hilbert
functions (and notice the shift, and that the lines for $R/C_1$
and $R/C_2$ are the first difference of those Hilbert functions,
i.e. the $h$-vectors):

\begin{center}
\begin{tabular}{c|cccccccccccccccccccccccc}
$R/I_s$ & & 1 & $s$ & \dots & $h_{R/I_s}\left ( \binom{s}{2} -1 \right )$ & $h_{R/I_s}\left ( \binom{s}{2} \right )$ & \dots \\
$R/C_1$ &  1 & $s-1$ & \dots & \dots & $\Delta h_{C_1}\left (
\binom{s}{2} \right )$ & $\Delta h_{C_1} \left ( \binom{s}{2} +1
\right )$ & \dots \\ \hline $R/G$ & 1 & $s$ & \dots & \dots & $A$
& $B$ & \dots
\end{tabular}
\end{center}
where
\begin{equation}\label{r/k}
\begin{array}{rcllll}
h_{R/G}\left ( \binom{s}{2} \right ) & = & A & = & h_{R/I_s} \left ( \binom{s}{2} -1 \right ) + \Delta h_{R/C_1} \left ( \binom{s}{2} \right ), \\  \\
h_{R/G}\left ( \binom{s}{2} +1 \right ) & = & B & = & h_{R/I_s}
\left ( \binom{s}{2}  \right ) + \Delta h_{R/C_1} \left (
\binom{s}{2} +1 \right ),
\end{array}
\end{equation}
and
\begin{center}
\begin{tabular}{c|cccccccccccccccccccccccc}
$R/G$ & & 1 & $s$ & \dots & $A$ & $B$ & \dots \\
$R/C_2$ &  1 & $s-1$ & \dots & \dots & $\Delta h_{C_2}\left ( \binom{s}{2} +1 \right )$ & $\Delta h_{C_2} \left ( \binom{s}{2} +2 \right )$ & \dots \\ \hline
$R/J_s$ & 1 & $s$ & \dots & \dots & $C$ & $D$ & \dots
\end{tabular}
\end{center}
where
\[
\begin{array}{rcllll}
h_{R/J_s}\left ( \binom{s}{2}+1 \right ) & = & C & = & h_{R/I_s} \left ( \binom{s}{2} -1 \right ) + \Delta h_{R/C_1} \left ( \binom{s}{2} \right ) + \Delta h_{R/C_2} \left ( \binom{s}{2} +1 \right ), \\  \\
h_{R/J_s}\left ( \binom{s}{2} +2 \right ) & = & D & = & h_{R/I_s} \left ( \binom{s}{2}  \right ) + \Delta h_{R/C_1} \left ( \binom{s}{2} +1 \right ) + \Delta h_{R/C_2} \left ( \binom{s}{2} +2 \right ).
\end{array}
\]
Now observe that the complete intersection $G$ has odd socle
degree $s(s-1)+1$; hence $A = B$.  Then it follows from
(\ref{r/k}) that
\begin{equation}\label{r/ir}
\hbox{$h_{R/I_s}\left ( \binom{s}{2} \right ) - h_{R/I_s} \left ( \binom{s}{2} -1 \right ) =
\Delta h_{R/C_1} \left ( \binom{s}{2} \right ) - \Delta h_{R/C_1} \left ( \binom{s}{2} +1 \right ).$}
\end{equation}
Thus we obtain
\begin{equation} \label{r/jr}
\begin{array}{rcl}
h_{R/J_s} (\binom{s}{2} +1) - h_{R/J_s} (\binom{s}{2}+2) & = &
h_{R/I_s}(\binom{s}{2} -1) - h_{R/I_s} (\binom{s}{2}) \\ \\
&& + \Delta h_{R/C_1}(\binom{s}{2}) - \Delta h_{R/C_1}(\binom{s}{2} +1) \\ \\
&& + \Delta h_{R/C_2}(\binom{s}{2} +1) - \Delta h_{R/C_2}(\binom{s}{2}+2) \\ \\
& = & \Delta h_{R/C_2}(\binom{s}{2} +1) - \Delta h_{R/C_2}(\binom{s}{2}+2).

\end{array}
\end{equation}
Combining (\ref{r/ir}) and (\ref{r/jr}), we see that it remains to show that
\begin{equation} \label{toshow}
\hbox{$\Delta h_{R/C_1} \left ( \binom{s}{2} \right ) - \Delta h_{R/C_1} \left ( \binom{s}{2} +1 \right ) \leq
\Delta h_{R/C_2}(\binom{s}{2} +1) - \Delta h_{R/C_2}(\binom{s}{2}+2). $}
\end{equation}
By the symmetry of the $h$-vectors of $R/C_1$ and $R/C_2$ we see that this is equivalent to showing
\begin{equation} \label{toshow2}
\hbox{$\Delta h_{R/C_1} \left ( \binom{s-1}{2} \right ) - \Delta h_{R/C_1} \left ( \binom{s-1}{2} -1 \right ) \leq
\Delta h_{R/C_2}(\binom{s-1}{2}) - \Delta h_{R/C_2}(\binom{s-1}{2}-1). $}
\end{equation}
Now, $\Delta h_{R/C_i}$ ($i = 1,2$) is the Hilbert function of an Artinian monomial complete intersection in $R$, namely $R/C_i'$, where $C_i'$ is obtained from $C_i$ by adding the missing variable.    Furthermore, if we replace $C_2'$ by $D_2 = (x_1,x_2^{s+1},x_3^s,\dots,x_s^s)$, we have that $R/C_2'$ and $R/D_2$ have the same Hilbert function, and $D_2 \subset C_1'$.

But such ideals have the Weak Lefschetz property (\cite{stanley}, \cite{watanabe}).  In particular, if $L$ is a general linear form, then the left-hand side of (\ref{toshow2}) is the Hilbert function of $R/(C_1' + (L))$  in degree $\binom{s-1}{2}$ and the right-hand side is the Hilbert function of $R/(D_2 + (L))$ in the same degree.  Because of the inclusion of the ideals, (\ref{toshow2}) follows and so the proof is complete.
\end{proof}

We now come to the main result of this section.  The case $r = 3$
was  proven by Brenner and Kaid \cite{BK}.  Note that when $r \leq
2$, all quotients of $R$ have the WLP by a result of \cite{HMNW}.

\begin{theorem} \label{Irr fails WLP}
Let $R = K[x_1,\dots,x_r]$, with $r \geq 3$, and consider
\[
I_{r,r} = (x_1^r,\dots,x_r^r,x_1x_2\cdots x_r ).
\]
Then the level Artinian algebra $R/I_{r,r}$ fails to have the WLP.
\end{theorem}

 \begin{proof}

 Specifically, we will check that the multiplication on $R/I_{r,r}$ by a general linear form fails
 surjectivity from degree $\binom{r}{2} - 1$ to degree $\binom{r}{2}$, even though the value of the Hilbert function is non-increasing between these two degrees.

 The proof is in two steps.

 \bigskip

 \noindent {\bf Step 1.}  We first prove this latter fact, namely that
 \[
 h_{R/I_{r,r}}(d-1) \geq h_{R/I_{r,r}}(d) \hbox{\ \ \  for $d = \binom{r}{2}$.  }
 \]

To do this, we again use basic double linkage.  We observe that
\[
\begin{array}{cclccccccccccccccc}
J_1 := (x_1^r, x_2^{r-1},\dots,x_r^{r-1}, x_1) & \leadsto  & x_2 \cdot (x_1^r, x_2^{r-1},\dots,x_r^{r-1},x_1) + (x_1^r, x_3^{r-1},\dots,x_r^{r-1}) \\ \\
& & = (x_1^r , x_2^r, x_3^{r-1} ,\dots, x_r^{r-1}, x_1 x_2) := J_2 \\ \\
& \leadsto & x_3 \cdot (x_1^r , x_2^r, x_3^{r-1} ,\dots, x_r^{r-1}, x_1 x_2) +
(x_1^r, x_2^r, x_4^{r-1}, \dots, x_r^{r-1}) \\ \\
&& = (x_1^r, x_2^r, x_3^r, x_4^{r-1} , \dots, x_r^{r-1}, x_1 x_2 x_3) := J_3 \\
& \vdots \\
& \leadsto & x_r \cdot (x_1^r, \dots, x_{r-1}^r, x_r^{r-1}, x_1 \cdots x_{r-1}) + (x_1^r, \dots, x_{r-1}^r) \\ \\
&& = (x_1^r, \dots, x_r^r, x_1 \cdots x_r) := J_r = I_{r,r}
\end{array}
\]
and we note that the first ideal, $J_1$,  is just $(x_2^{r-1},\dots,x_r^{r-1},x_1)$.  Furthermore, this ideal is a complete intersection with socle degree $(r-1)(r-2)$.  The midpoint of the Hilbert function is in degree $\binom{r-1}{2}$.

Let $C_1 = (x_1^r, x_3^{r-1}, \dots, x_r^{r-1})$ and $C_2 = (x_1^r, x_2^r, x_4^{r-1}, \dots, x_r^{r-1})$.  We note that the first difference of the Hilbert function of $R/C_1$ is symmetric with odd socle degree,  so the values in degrees $\binom{r-1}{2}$ and $\binom{r-1}{2} +1$ are equal.  The key point, though, is that by applying Lemma \ref{diff of hf} with $s = r-1$, we obtain that
\begin{equation} \label{keypt}
\hbox{$ h_{R/J_1}(\binom{r-1}{2}) - h_{R/J_1} ( \binom{r-1}{2} -1) \leq h_{R/C_2}(\binom{r-1}{2} +1) - h_{R/C_2}(\binom{r-1}{2} +2).$}
\end{equation}

We are interested in the values of the Hilbert function of $R/I_{r,r}$ in degrees $\binom{r}{2} -1$ and $\binom{r}{2}$.  Since $I_{r,r}$ is obtained from $J_1$ by a sequence of $r-1$ basic double links, and since each one involves a shift by 1 of the Hilbert function, this corresponds (first) to an examination of the Hilbert function of $R/J_1$ in degrees $\binom{r}{2}-1 - (r-1) = \binom{r-1}{2} - 1$ and $\binom{r-1}{2}$.  The former is smaller than the latter, but we do not need to know the precise values.

Our observation  in the paragraph preceding   (\ref{keypt})  shows that when we add the first difference of the (shifted) Hilbert function of $R/C_1$ to get $h_{R/J_2}(\binom{r-1}{2})$ and $h_{R/J_2}(\binom{r-1}{2} +1)$, the difference between these two values is the same as the difference between the values of the Hilbert function of $R/J_1$ in degrees $\binom{r-1}{2}-1$ and $\binom{r-1}{2}$, with the latter being larger.  However, the point of (\ref{keypt}) is that when we then add the (shifted) first difference of the Hilbert function of $R/C_2$, we overcome this difference and already have a Hilbert function with the value in degree $\binom{r-1}{2} + 1$ larger than that in degree $\binom{r-1}{2} +2$.  Since each subsequent Hilbert function has the (shifted) value in the smaller of the corresponding degrees larger than the value in the second, we finally obtain the same for the desired Hilbert function, namely that of $R/I_{r,r}$.  This concludes step 1.

 \vskip 2mm

 \noindent {\bf Step 2.}  To prove that $R/I_{r,r}$ fails surjectivity from
degree $\binom{r}{2} - 1$ to degree $\binom{r}{2}$, we will use Proposition \ref{splitting and surj} (d).  Note that
\[
\bar I_{r,r} \cong (x_1^r,\dots,x_{r-1}^r, (x_1 + \dots + x_{r-1})^r, x_1 \cdots x_{r-1} \cdot (x_1 + \cdots + x_{r-1}) ).
\]
  We now claim that it is enough to verify that there is a homogeneous form $F\in \bar R \cong K[x_1, \cdots ,x_{r-1}]$ of degree $\binom{r-1}{2}$ such that

\begin{equation}\label{First equation}
F \cdot x_1x_2 \cdots x_{r-1}(x_1+ \cdots +x_{r-1}) \in (x_1^r, \dots ,x_{r-1}^r)
\end{equation}
and
\begin{equation}\label{Second equation}
F \cdot (x_1 + ... + x_{r-1})^r \in (x_1^r,...,x_{r-1}^r).
\end{equation}

Indeed, clearly $x_1^r ,\dots , x_{r-1}^r$ is a regular sequence
in $\bar I_{r,r}$ that extends to a minimal generating set for
$I_{r,r}$. So Proposition \ref{splitting and surj} (d) shows that
it is enough to find a form $F$ of degree $r(r-1) - \left (
\binom{r}{2} + r-1 \right ) = \binom{r-1}{2}$, non-zero modulo
$(x_1^r,\dots,x_{r-1}^r)$, such that (\ref{First equation}) and
(\ref{Second equation}) hold.  Hence our claim holds.

The heart of the proof is to show that the specific polynomial
\[
F=\sum _{i_1+\cdots +i_{r-1} = \binom{r-1}{2} \atop
0\le i_j \le r-2 , \quad i_j\ne i_{\ell }} (-1)^{sgn(i_1,\cdots
,i_{r-1})}x_1^{i_1}\cdots x_{r-1}^{i_{r-1}}
\]
satisfies (\ref{First equation}) and (\ref{Second equation}).  Note that $F$ is the determinant of a Vandermonde matrix.   $F$~simply  consists of a sum of terms, all with coefficient 1 or $-1$, obtained as follows.  Each term consists of a product of different powers of the $r-1$ variables (remember that we are in the quotient ring).  Namely, for {\em each} permutation, $\sigma$, of $(0,1,2,...,r-2)$, we look at the term
\[
      (-1)^{sgn(\sigma)} \cdot A
\]
where A is the monomial obtained by taking the $i$-th variable to the power given by the $i$-th entry in $\sigma$.  For example, if $r = 5$ and $\sigma = (2,0,3,1)$ then $sgn(\sigma) = -1$ so we have the term
$ - x_1^2 x_3^3 x_4$.

We first check (\ref{First equation}). In order for the product to
be contained in $(x_1^r,\dots,x_{r-1}^r)$, we need that every term
in the product that does not contain at least one
  exponent $\geq r$ be canceled by another term in the product.  That is,
we have
 $F \cdot x_1x_2...x_{r-1}(x_1+...+x_{r-1}) \in (x_1^r,...,x_{r-1}^r)$ if and only if
\[
\sum _{m = 1}^{r-1} \left (
 \sum _{i_1+\cdots +i_{r-1}={r-1\choose 2} \atop
0\le i_j \le r-2, \quad i_j\ne i_{\ell }, \quad i_{m}\ne r-2}
(-1)^{sgn(i_1,\cdots , i_{r-1})}x_1^{i_1+1} \cdots ,
x_{m}^{i_{m}+2}, \cdots ,i_{r-1}^{i_{r-1}+1} \right )=0.
\]

\vskip 4mm Notice that we have ruled out $i_m = r-2$ since otherwise $i_m + 2 = r$ and that term is automatically in the desired ideal.  But then the hypotheses $i_1+ \cdots
+i_{r-1}={r-1\choose 2}$, $0\le i_j \le r-2$ and $i_j \ne i_{\ell}$ imply that there exists a unique  integer $n$ with $1 \le n \le r-1$ such that $i_{m}+2=i_{n}+1$. Hence the summand
\[
(-1)^{sgn(i_1,\cdots , m , \cdots , n , \cdots , i_{r-1})}x_1^{i_1+1} \cdots ,
x_{m}^{i_{m}+2}, \cdots ,x_{r-1}^{i_{r-1}+1}
\]
 is cancelled against the summand
 \[
 (-1)^{sgn(i_1,\cdots , m+1, \cdots , n - 1,
\cdots , i_{r-1})}x_1^{i_1+1} \cdots , x_{m}^{i_{m}+2}, \cdots
,x_{r-1}^{i_{r-1}+1}.
 \]
(For notational convenience we have assumed $m < n$ but this is not at all important.)

\vskip 4mm We now  prove (\ref{Second equation}). We have
 \[
  F \cdot (x_1 + ... + x_{r-1})^r = F\cdot \left ( \sum_{j_1+\cdots
+j_{r-1}=r \atop j_i\ge 0 }\frac{r!}{j_1!\cdots
j_{r-1}!}x_1^{j_1}\cdots x_{r-1}^{j_{r-1}} \right ) =
  \]
  \[
\left ( \sum_{j_1+\cdots +j_{r-1}=r \atop j_i\ge 0 }\frac{r!}{j_1!\cdots
j_{r-1}!}x_1^{j_1}\cdots x_{r-1}^{j_{r-1}} \right ) \cdot \left ( \sum
_{i_1+\cdots +i_{r-1}={r-1\choose 2} \atop 0\le i_j \le r-2 ,
\quad i_j\ne i_{\ell }} (-1)^{sgn(i_1,\cdots
,i_{r-1})}x_1^{i_1}\cdots x_{r-1}^{i_{r-1}} \right ) =
  \]

\[
  \sum_{j_1+\cdots +j_{r-1}=r \atop j_i\ge 0} \left ( \sum _{i_1+\cdots
+i_{r-1}={r-1\choose 2} \atop 0\le i_j \le r-2 , \quad i_j\ne
i_{\ell }}\frac{(-1)^{sgn(i_1,\cdots ,i_{r-1})}r!}{j_1!\cdots
j_{r-1}!}x_1^{i_1+j_1}\cdots x_{r-1}^{i_{r-1}+j_{r-1}} \right ).
  \]
Therefore
  \[
   F \cdot (x_1 + ... + x_{r-1})^r \in (x_1^r,...,x_{r-1}^r)
   \]
  if and only if
  \[
   G := \sum_{j_1+\cdots +j_{r-1}=r \atop r-1\ge j_i\ge 0} \left ( \sum _{i_1+\cdots
+i_{r-1}={r-1\choose 2} \atop 0\le i_{\ell } \le
\min(r-2,r-1-j_{\ell}) , \quad i_j\ne i_{\ell
}}\frac{(-1)^{sgn(i_1,\cdots ,i_{r-1})}r!}{j_1!\cdots
j_{r-1}!}x_1^{i_1+j_1}\cdots x_{r-1}^{i_{r-1}+j_{r-1}} \right )=0.
   \]

\vskip 2mm Given an $(r-1)$-uple of non-negative integers
$\underline{j}:=(j_1,\cdots ,j_{r-1})$ such that $j_1+\cdots
+j_{r-1}=r$, we set
\[
   C_{\underline{j}}:= \frac{r!}{j_1!\cdots j_{r-1}!}.
\]
Notice that two $(r-1)$-uples of non-negative integers $(j_1,\cdots , j_{r-1})$ and $(j'_1,\cdots ,j'_{r-1})$ with $j_1 + \cdots +j_{r-1} = r = j'_1 + \cdots + j'_{r-1}$ verify
\begin{equation}
\frac{r!}{j_1! \cdots j_{r-1}!} = \frac{r!}{j'_1!\cdots j'_{r-1}!} \Leftrightarrow \{j_1, \cdots , j_{r-1} \}  = \{j'_1,\cdots, j'_{r-1} \}.
\end{equation}

\vskip 2mm Given an $(r-1)$-uple of non-negative integers
$\underline{j} := (j_1,\cdots ,j_{r-1})$ such that $j_1+\cdots + j_{r-1} = r$ and an $(r-1)$-uple
$\underline{\alpha  }:=(\alpha_1, \cdots , \alpha _{r-1})$, we define (from now on $\# (B)$
means the cardinality of the set $B$):
\[
N_{\underline{\alpha }, \underline{j}}^{} : = \# (A(\underline{\alpha})_{\underline{j}})
\]
 where $A(\underline{\alpha})_{\underline{j}}$ is the set of monomials
$\pm C_{\underline{j}} x_1^{\alpha_1} \cdots x_{r-1}^{\alpha_{r-1}}$ in $G$ of multidegree  $\underline{\alpha}$ and coefficient $\pm C_{\underline{j}}$. To prove (\ref{Second
equation}), it is enough to see that $N_{\underline{\alpha},\underline{j}}$ is even and half of the elements of $A(\underline{\alpha })_{\underline{j}}$ have coefficient $+C_{\underline{j}}$ and the other half have coefficient $-C_{\underline{j}}$.  Let us prove it. Without loss of generality
we can assume that $\alpha _1 \ge \alpha _2 \ge \cdots \ge \alpha_{r-1}$ (we re-order the variables, if necessary). We will see that for any monomial in $A(\underline{\alpha})_{\underline{j}}$
there is a well determined way to associate another monomial in
$A(\underline{\alpha})_{\underline{j}}$ with the opposite sign.  Indeed, the monomials in $A(\underline{\alpha})_{\underline{j}}$ have degree ${r-1\choose 2}+r={r\choose 2} +1$ and, moreover, $0 \le
\alpha_{\ell} \le r-1$ for all $1\le \ell \le r-1$.
Therefore, there exist integers $1\le p<q\le r-1$ such that $\alpha_p = \alpha_q$.

We define $p_0 := \min\{ p \mid \alpha_{p} = \alpha_{p+1} \}$. Now, we take an arbitrary monomial
\[
 \pm C_{\underline{j}} x_1^{\alpha_1} \cdots x_{p_0}^{\alpha_{p_{0}}} x_{p_0 +1}^{\alpha_{p_0+1}}\cdots  x_{r-1}^{\alpha _{r-1}}\in A(\underline{\alpha})_{\underline{j}}
 \]
  where $\alpha_1=j_1+i_1$, $\dots$, $\alpha_{p_{0}} = j_{p_0} + i_{p_0}$, $\alpha_{p_0+1} = j_{p_0+1} + i_{p_0 +1}$, $\dots$, $\alpha_{r-1} = j_{r-1} + i_{r-1}$. It will be cancelled with
  \[
  \mp C_{\underline{j}} x_1^{\alpha_1} \cdots x_{p_0}^{\alpha_{p_{0}}} x_{p_0 +1}^{\alpha_{p_0+1}}\cdots  x_{r-1}^{\alpha _{r-1}}\in A(\underline{\alpha})_{\underline{j}}
  \]
 where $\alpha_1=j_1+i_1$, $\dots$ ,$\alpha_{p_{0}} = j_{p_0+1} + i_{p_0 +1}$, $\alpha_{p_0+1} = j_{p_0} + i_{p_0}$, $\dots$ , $\alpha_{r-1} = j_{r-1} + i_{r-1}$ and we are done.
\end{proof}

\begin{example}
We illustrate the construction in Step 1 of the proof of Theorem
\ref{Irr fails WLP} for the case $r=5$.  In the following table of
Hilbert functions and $h$-vectors, we have
\[
\begin{array}{rcl}
J_1 & =  & (x_1,x_2^4,x_3^4,x_4^4,x_5^4) \\
J_2 & = & x_2 \cdot J_1 + (x_1^5,x_3^4,x_4^4,x_5^4) = (x_1^5,x_2^5,x_3^4,x_4^4,x_5^4,x_1x_2) \\
J_3 & = & x_3 \cdot J_2 + (x_1^5,x_2^5,x_4^4,x_5^4) = (x_1^5,x_2^5,x_3^5,x_4^4,x_5^4,x_1x_2x_3) \\
J_4 & = & x_4 \cdot J_3 + (x_1^5,x_2^5,x_3^5,x_5^4) = (x_1^5,x_2^5,x_3^5, x_4^5,x_5^4,x_1x_2x_3x_4) \\
J_5 & = & x_5 \cdot J_4 + (x_1^5,x_2^5,x_3^5,x_4^5) = (x_1^5,x_2^5,x_3^5,x_4^5,x_5^5,x_1x_2x_3x_4x_5)
\end{array}
\]

In the following calculation, we have put in boldface the critical range of degrees.

{\small

\begin{tabular}{c|ccccccccccccccccccccccccccc}
Ideal & \multicolumn{17}{c}{Hilbert function/$h$-vector}
 \\ \hline
$J_1$ & && & &   1 & 4 & 10 & 20 & 31 & {\bf 40} & {\bf  44} & 40 & 31 & 20 & 10 & 4 & 1  \\
(5,4,4,4) &&&& 1 & 4 & 10 & 20 & 32 & 43 & {\bf 50} & {\bf 50} & 43 & 32 & 20 & 10 & 4 & 1 \\ \hline
$J_2$ &&&& 1 & 5 & 14 & 30 & 52 & 74 & {\bf 90} & {\bf 94} & 83 & 63 & 40 & 20 & 8 & 2 \\
(5,5,4,4) &&& 1 & 4 & 10 & 20 & 33 & 46 & 56 & {\bf 60} & {\bf 56} & 46 & 33 & 20 & 10 & 4 & 1 \\ \hline
$J_3$ &&& 1 & 5 & 15 & 34 & 63 & 98 & 130 & {\bf 150} & {\bf 150} & 129 & 96 & 60 & 30 & 12 & 3 \\
(5,5,5,4) && 1& 4& 10& 20& 34& 49& 62& 70& {\bf 70}& {\bf 62} & 49& 34& 20& 10& 4& 1 \\ \hline
$J_4$ && 1& 5& 15& 35& 68& 112& 160& 200& {\bf 220} & {\bf 212} & 178& 130& 80& 40& 16& 4 \\
(5,5,5,5) & 1& 4& 10& 20& 35& 52& 68& 80& 85& {\bf 80} & {\bf 68} & 52& 35& 20& 10& 4& 1 \\ \hline
$J_5$ & 1& 5& 15& 35& 70& 120& 180& 240& 285& {\bf 300} & {\bf 280} & 230& 165& 100& 50& 20& 5
\end{tabular}
}

\bigskip

It is interesting to note that experimentally we have verified
that $R/J_1$ and $R/J_2$ have the WLP, while $R/J_3$, $R/J_4$ and
$R/J_5$ do not.  The algebras that fail to have the WLP all fail
surjectivity in the range indicated in boldface.  Only $R/J_5$
fails to have the WLP in any other degree, namely it fails
injectivity in the preceding degree.
\end{example}

As mentioned above, we now have a partial answer to Question 4.2 of \cite{MMR3}.  Recall that $A(n)$ is defined to be the minimum number (if it exists) such that {\em every} Artinian ideal $I \subset K[x_1,\dots,x_n]$ with number of generators $\mu(I) \leq A(n)$ has the WLP.

\begin{corollary}
 If $A(n)$ exists then it equals $n$.
\end{corollary}


\section{An almost monomial almost complete intersection} \label{almost monomial}

In order to illustrate the subtlety of the Weak Lefschetz Property, we now describe a class of ideals that is very similar to the class of ideals discussed in Section \ref{arb codim}.  That is, we consider, for each codimension $r \geq 3$, the ideal
\[
\mathfrak{J}_r = (x_1^r,\dots,x_r^r,x_1\cdots x_{r-1}(x_1 + x_r)).
\]
We will compare the properties of this ideal with those of the ideal
\[
I_{r,r} = (x_1^r,\dots,x_r^r,x_1\cdots x_r).
\]
Included in this subtlety is the fact that  the WLP behavior
changes with the characteristic. Notice that our results in
positive characteristic do not depend on whether the field is
finite or not.

Our first result shows that we cannot distinguish the two ideals by
solely looking at their Hilbert functions.

\begin{lemma}
$R/\mathfrak{J}_r$ and $R/I_{r,r}$ have the same Hilbert function.
\end{lemma}

\begin{proof}
We will show that $\mathfrak{J}_r$ arises via a sequence of basic double links which are numerically equivalent to the one that produced $I_{r,r}$ in the first part of
Theorem \ref{Irr fails WLP}.  Notice first that the ideals
\[
(x_1^r, x_2^{r-1},\dots,x_r^{r-1}, x_1) = (x_2^{r-1},\dots,x_r^{r-1}, x_1) \ \ \ \hbox{and} \ \ \ (x_1^{r-1},\dots,x_{r-1}^{r-1}, x_1+x_r)
\]
have the same Hilbert function.
Notice also that this latter ideal is equal to
\[
(x_1^{r-1},\dots,x_{r-1}^{r-1}, x_r^r, x_1+x_r).
\]
 In Step 1 of Theorem \ref{Irr fails WLP} we saw a sequence of basic double links starting with the ideal $(x_1^{r}, x_2^{r-1},\dots,x_{r}^{r-1}, x_1)$ and ending with $(x_1^r,\dots,x_r^r,x_1,\dots,x_r)$.  We will now produce a parallel sequence of basic double links starting with $(x_1^{r-1},\dots,x_{r-1}^{r-1}, x_r^r, x_1+x_r)$ and ending with $(x_1^r,\dots,x_r^r,x_1\cdots x_{r-1}(x_1+x_r))$, such that at each step the two sequences are numerically the same, and hence the resulting ideals at each step have the same Hilbert function.
\[
\begin{array}{lcclccccccccccccccc}
(x_1^{r-1}, \dots,x_{r-1}^{r-1},x_r^r, x_1+x_r) \\ \\
\begin{array}{ll}
 \leadsto  & x_{r-1}
\cdot (x_1^{r-1}, \dots,x_{r-1}^{r-1},x_r^r, x_1+x_r)
+ (x_1^{r-1}, \dots, x_{r-2}^{r-1}, x_r^{r}) \\ \\

 & = (x_1^{r-1} , x_2^{r-1} ,\dots, x_{r-2}^{r-1}, x_{r-1}^r, x_r^{r}, (x_1+x_r)\cdot x_{r-1})  \\ \\

 \leadsto &  x_{r-2} \cdot (x_1^{r-1} , x_2^{r-1} ,\dots, x_{r-2}^{r-1}, x_{r-1}^r, x_r^{r}, (x_1+x_r)\cdot x_{r-1})    + (x_1^{r-1}, \dots, x_{r-3}^{r-1}, x_{r-1}^r, x_r^{r}) \\ \\

& = (x_1^{r-1} , x_2^{r-1} ,\dots, x_{r-3}^{r-1},x_{r-2}^r,  x_{r-1}^r, x_r^{r}, (x_1+x_r)\cdot x_{r-2} x_{r-1})  := J_3 \\
&  \vdots \\  \\

 \leadsto &  x_1 \cdot (x_1^{r-1}, x_2^r, \dots, x_{r-1}^r, x_r^{r}, (x_1+x_r)\cdot x_2 \cdots x_{r-1}) + (x_2^r, \dots, x_{r}^r) \\ \\
 & = (x_1^r, \dots, x_r^r, x_1 \cdots x_{r-1}(x_1+x_r)) = {\mathfrak
 J}_r.
 \end{array}
\end{array}
\]
This completes the proof.
\end{proof}

We will now show that the two algebras behave differently with respect to the WLP. Recall that $R/I_{r, r}$ does not have the WLP if $r \geq 3$. Studying $R/\mathfrak{J}_r$ when $r = 3$ is not too difficult:

\begin{proposition}
\label{prop:r3-aci}
For every field $K$, the algebra
\[
R/\mathfrak{J}_3 = K[x, y, z]/(x^3, y^3, z^3, xy(x+z))
\]
 has the WLP if and only if  the characteristic of $K$ is not three.
\end{proposition}

\begin{proof}
If the characteristic of $K$ is three then for  every linear form
$\ell \in R$, $\ell^3$ is in $(x^3, y^3, z^3)$.   Thus the residue
class of $\ell^2$ is in the kernel of the multiplication map
\[
\times \ell: (R/\mathfrak{J}_3)_2 \to (R/\mathfrak{J}_3)_3.
\]
This shows that $R/\mathfrak{J}_3$ does not have  the WLP if $\chara
K = 3$.

Now assume that $\chara K \neq 3$. Consider the linear form $L =
x+y+z$. Then one checks that
\[
(\mathfrak{J}_3, L)/(L) \cong ((x, y)^3, L)/(L),
\]
which implies that the multiplication map
\[
\times L: (R/\mathfrak{J}_3)_2 \to (R/\mathfrak{J}_3)_3
\]
is surjective. Hence $R/\mathfrak{J}_3$ has the WLP in this case.
\end{proof}

The case when $r = 4$ is considerably more complicated.

\begin{proposition}
\label{prop:r4-aci} For every field $K$, the algebra
\[
R/\mathfrak{J}_4 = K[w,x, y, z]/(w^4, x^4, y^4, z^4, wxy(w+z))
\]
 has the WLP if and only if  the characteristic of $K$ is not two or five.
\end{proposition}

\begin{proof}
The Hilbert function of $R/\mathfrak{J}_4$ is $1, 4, 10, 20, 30, 36, 34, \ldots$.
Hence, by Proposition \ref{gen wlp},  $R/\mathfrak{J}_4$ has the WLP if and only if, for a general form $L$,
the multiplication maps
\[
\times L: (R/\mathfrak{J}_4)_4 \to (R/\mathfrak{J}_4)_5 \quad
\text{and} \quad \times L: (R/\mathfrak{J}_4)_5 \to
(R/\mathfrak{J}_4)_6
\]
are injective and surjective, respectively.

We first show that the latter map is surjective if $L := 2 w + x +
y + z$, provided the characteristic of $K$ is neither 2 nor 5.
Notice that this map is surjective if and only if
$(R/(\mathfrak{J}_4, L)])_6 = 0$. Since
\[
R/(\mathfrak{J}_4, L) \cong K[w, x, y]/(w^4, x^4, y^4, (2w + x + y)^4, wxy(w+x+y))
\]
this is equivalent to the fact that $\dim_K ((w^4, x^4, y^4, (2w +
x + y)^4, wxy(w+x+y))_6 = 28$. To compute the dimension of $((w^4,
x^4, y^4, (2w + x + y)^4, wxy(w+x+y))_6$, we consider the
coefficients of the 28 degree 6 monomials in $K[w, x, y]$
occurring in each of the  30 polynomials $f q$, where $f$ is one
of the forms $w^4, x^4, y^4, (2w + x + y)^4, wxy(w+x+y)$ and $q$
is one of the quadrics $w^2, wx, x^2, wy, xy, y^2$. Compute these
coefficients assuming, temporarily, that $\chara K = 0$, and
record them in a $30 \times 28$ matrix $M$ whose entries are
integers.  Using CoCoA we verified that the greatest common
divisor of all the maximal minors of $M$ is $320 = 2^8 \cdot 5$.
This shows that the matrix $M$ has rank 28 if and only if $\chara
K \neq 2, 5$.

We now discuss the map $\times L:  (R/\mathfrak{J}_4)_4 \to
(R/\mathfrak{J}_4)_5$, where $L$ is a general linear form. This
map is injective if and only if $\dim_K (R/(\mathfrak{J}_4, L))_5
= 6$,
which is equivalent to \\
$\dim_K ((\mathfrak{J}_4, L)/(L))_5 = 15$.

Assume first that the field $K$ is infinite. Then an argument
similar to the one in the proof of Proposition \ref{lem-L-element}
shows we may assume that
\[
L := t w +  x+ y -  z,
\]
where $t \in K$. Then
\[
R/(\mathfrak{J}_4, L) \cong
K[w, x, y]/(w^4, x^4, y^4, (t w + x + y)^4, wxy((t+1) \cdot w+x+y)).
\]
To compute the dimension of $((w^4, x^4, y^4, (t w + x + y)^4,
wxy((t+1) \cdot w+x+y)))_5$, we consider the coefficients of the
21 degree 5 monomials in $K[w, x, y]$ occurring in each of the  15
polynomials $f \ell$, where $f$ is one of the forms $w^4, x^4,
y^4, (t w + x + y)^4, wxy((t+1) \cdot w+x+y)$ and $\ell$ is one of
the variables $w, x, y$. Compute these coefficients assuming,
temporarily, that $\chara K = 0$, and record them in a $15 \times
21$ matrix $N$ whose entries are polynomials in $\Z[t]$. A CoCoA
computation provides that all maximal minors of $N$ are divisible
by 10 and that one of the minors is $80 t^4 (t+1)^2$. It follows
that the rank of $N$ is 15 if and only if $\chara K \neq 2, 5$.
Hence we have shown that over an infinite field, for a general
linear form $L$, the map $\times L: (R/\mathfrak{J}_4)_4 \to
(R/\mathfrak{J}_4)_5$ is injective if and only if the
characteristic of $K$ is neither 2 nor 5. This also implies that
$R/\mathfrak{J}_4$ does not have the WLP if $K$ is a finite field
of characteristic 2 or 5. Furthermore, every field whose
characteristic is not 2 or 5 contains an element $t$ such that $80
t^4 (t+1)^2$ is not zero. Hence the above arguments show that in
this case there is a linear form $L$ such that $\times L:
(R/\mathfrak{J}_4)_4 \to (R/\mathfrak{J}_4)_5$ is injective.

Combining this with the first part of the proof, our assertion
follows.
\end{proof}

\begin{remark}
  \label{rem:conj-wlp}

(i) For $R/\mathfrak{J}_4$, Proposition \ref{prop:r4-aci} provides
an affirmative answer to Problem \ref{prop-ext}.

(ii) We expect that in characteristic zero, for each integer  $r
\geq 2$, the algebra $R/\mathfrak{J}_r$ has the WLP and that $L = 2
x_1 + x_2 + \cdots + x_{r-1} - x_r$ is a Lefschetz element.
\end{remark}


\section{Monomial almost complete intersections in three variables}
\label{monomial aci}

Now we consider ideals of the form
\[
I = (x^a, y^b, z^c, x^\alpha y^\beta z^\gamma )
\]
in $R = K[x,y,z]$, where $0 \leq \alpha < a$, $0 \leq \beta < b$ and $0 \leq \gamma < c$.
This class of ideals was first considered in \cite{brenner}, Corollary 7.3.

\begin{proposition} \label{codim3}
If $I$ is as above and is not a complete intersection  then

\begin{itemize}
\item[(i)] The inverse system for $I$ is given by $(x^{a-1} y^{b-1} z^{\gamma -1}, x^{a-1} y^{\beta-1} z^{c-1}, x^{\alpha -1} y^{b-1} z^{c-1})$, where we make the convention that if a term has an exponent of $-1$ (e.g. if $\gamma = 0$), that term is removed.

\item[(ii)] In particular, if $\alpha, \beta, \gamma > 0$ then $R/I$ has Cohen-Macaulay type 3.  Otherwise it has Cohen-Macaulay type 2.

\item[(iii)] If $\alpha, \beta, \gamma > 0$, the socle degrees of $R/I$ are $b+c+\alpha -3, a+c+\beta-3, a+b+\gamma-3$.  In particular, $R/I$ is level if and only if $a-\alpha = b - \beta = c - \gamma$.

\item[(iv)] Suppose that one of $\alpha, \beta, \gamma = 0$; without loss of generality say $\gamma = 0$. Then the corresponding socle degree in (iii), namely $a+b+\gamma-3$, does not occur.  Now $R/I$ is level if and only if $a - \alpha = b - \beta$, and $c$ is arbitrary.

\item[(v)] Let $J = I : x^\alpha y^\beta z^\gamma$ be the ideal residual to $I$ in the complete intersection $(x^a, y^b, z^c)$.  Then $J = (x^{a-\alpha}, y^{b-\beta}, z^{c-\gamma})  $.

\item[(vi)] A free resolution of $R/I$ is
\begin{equation} \label{res of aci}
{\footnotesize 0 \rightarrow
\begin{array}{c}
R(-\alpha-b-c) \\
\oplus \\
R(-\beta -a-c) \\
\oplus \\
R(-\gamma -a-b)
\end{array}
\rightarrow
\begin{array}{c}
R(-\alpha -\beta - c) \\
\oplus \\
R(-\alpha - \gamma -b) \\
\oplus \\
R(-\beta - \gamma -a) \\
\oplus \\
R(-b-c) \\
\oplus \\
R(-a-c) \\
\oplus \\
R(-a-b)
\end{array}
\rightarrow
\begin{array}{c}
R(-\alpha - \beta - \gamma) \\
\oplus \\
R(-a) \\
\oplus \\
R(-b) \\
\oplus \\
R(-c)
\end{array}
 \rightarrow R \rightarrow R/I \rightarrow 0.}
\end{equation}
This is minimal if and only if $\alpha, \beta, \gamma $ are all positive.
\end{itemize}

\end{proposition}

\begin{proof}
Part (i) follows by inspection.  Then (ii), (iii) and (iv) follow
immediately from (i).  As before, (v) is a simple computation of the
colon ideal, based on the fact \cite{beintema} that $J$ is a
complete intersection, so it only remains to check the degrees.
Having (v), it is a straightforward computation using the mapping
cone to obtain (vi).
\end{proof}

\begin{theorem} \label{zero mod 3}
Assume that $K = \overline{K}$ is an algebraically  closed field
of characteristic zero.  For $I = (x^a, y^b, z^c, x^\alpha y^\beta
z^\gamma )$, if the WLP fails then $a+b+c+\alpha + \beta + \gamma
\equiv 0$ (mod 3).
\end{theorem}

\begin{proof}
Let $\mathcal E$ be the syzygy bundle  of $I$ and let $L \cong
\mathbb P^1$ be a general line.    By \cite{BK} Theorem 3.3, if
the WLP fails then $\mathcal E$ is semistable.  Furthermore,  the
splitting type of ${\mathcal E}_{\rm norm}$ must be $(1,0,-1)$
(apply \cite{BK}, Theorem  2.2 and the Grauert-M\"ulich theorem).
Hence the twists of ${\mathcal E}|_L$ are three consecutive
integers.  Since the restriction of $\mathcal E$ to $L$ is the
(free) syzygy module corresponding to the restriction of the
generators of $I$ to $\mathbb P^1$, we see that the sum of the
generators must be divisible by 3.
\end{proof}

\begin{corollary} \label{mod 3}
Assume that $K = \overline{K}$ is an algebraically closed field of
characteristic zero.  If $R/I$ is level and the WLP fails then $a
+ b + c \equiv 0$ (mod 3) and $\alpha + \beta + \gamma \equiv 0$
(mod 3).
\end{corollary}

\begin{proof}
Since $R/I$ is level, by Proposition \ref{codim3} we can write $a = \alpha +t$, $b = \beta + t$, $c = \gamma + t$ for some $t \geq 1$.  By Theorem \ref{zero mod 3}, we have
\[
2(\alpha + \beta + \gamma) + 3t \equiv 0 \hbox{ (mod 3)}.
\]
It follows that $\alpha + \beta + \gamma \equiv 0$ (mod 3), so again by Theorem \ref{zero mod 3} we also get $a+b+c \equiv 0$ (mod 3).
\end{proof}

\begin{remark}
The proof of Theorem \ref{zero mod 3} applies not only to monomial
ideals. Indeed, for any almost complete intersection in $R =
K[x,y,z]$, if the WLP fails then $\sum_{i=1}^4 d_i \equiv 0$ (mod
3), where $d_1,d_2,d_3,d_4$ are the degrees of the minimal
generators.
\end{remark}

A very interesting class of ideals is the following, recalling the notation introduced in \eqref{eq:def-ideal}.

\begin{corollary} \label{1st WLP}
The algebra $A = R/I_{3,k,3}$ has the following properties.

\begin{itemize}
\item[(a)] the socle degree is $e = 2k-2$.

\item[(b)] the peak of the Hilbert function occurs in degrees $k-1$ and $k$, and has value $3k-3$.

\item[(c)]  The corresponding inverse system is $(x^{k-1}y^{k-1},x^{k-1}z^{k-1},y^{k-1}z^{k-1})$.

\item[(d)] Assuming $\chara K \neq 2$, the WLP fails if and only
if $k$ is odd.  Note that in this case $e \equiv 0$ (mod 4).

\end{itemize}

\end{corollary}

\begin{proof}
Parts (a), (b) and (c) are immediate from Proposition \ref{codim3}. Part (d) is a special case of Theorem \ref{thm-wlp}.
\end{proof}

For the remainder of this section we focus on the WLP.  If the
non-pure monomial involves only two of the variables, then  the
ideal always has the WLP.

\begin{lemma}
  \label{lem:variable-missing}
Adopt the above notation. If $\al = 0$ and $K$ has characteristic
zero, then $R/I$ has the WLP.
\end{lemma}

\begin{proof}
The assumption $\al = 0$ provides that $R/I$ is isomorphic to $B
  \otimes C$, where $B = K[y, z]/(y^b, y^{\be} z^{\ga}, z^c)$ and $C
  =  K[x]/(x^a)$. By Proposition 4.4 in \cite{HMNW}, $B$ and $C$ have
the WLP, hence $A$ has the WLP by \cite{watanabe}, Corollary 3.5.
\end{proof}

If the non-pure monomial involves all three variables then, due to
Theorem \ref{codim3}(iii), the ideal is level if and only if $I$ is
of the form
\begin{equation}
  \label{eq-level-aci}
I_{\alpha, \beta, \gamma, t} =
(x^{\alpha + t} , y^{\beta + t}, z^{\gamma + t}, x^\alpha y^\beta z^\gamma),
\end{equation}
where $t \geq 1$ and, without loss of generality, $1 \leq \al \leq \be \leq \ga$.

Next, we analyze  when the syzygy bundle of $I$ is semistable.
(The relevance of semistability to the WLP was introduced in
\cite{HMNW} and generalized in \cite{BK}.)

\begin{lemma}
  \label{lem-sst}
Assume that $K$ is algebraically closed of characteristic zero.
Then the syzygy bundle of $I$ is semistable if and only if $\ga \leq 2 (\al + \be)$ and $\frac{1}{3} (\al + \be + \ga) \leq t$.
\end{lemma}

\begin{proof}
Brenner (\cite{brenner}, Corollary 7.3) shows that in general for an ideal $I =  (x^a, y^b, z^c, x^\alpha y^\beta z^\gamma )$, the syzygy bundle is semistable if and only if
\begin{itemize}
\item[(i)] $ 3 \max \{ a,b,c,\alpha + \beta + \gamma \} \leq a + b + c + \alpha + \beta + \gamma $, and

\item[(ii)] $\displaystyle \min \{ \alpha + \beta + c, \alpha + b + \gamma, a + \beta + \gamma, a+b, a+c, b+a \} \geq \frac{a+b+c+\alpha + \beta + \gamma}{3}$.
\end{itemize}
Applying this to our ideal $I$ condition (i) reads as
\[
3 \max \{\ga + t, \al + \be + \ga \} \leq 2 (\al + \be + \ga) + 3 t.
\]
Hence, it is equivalent to
\begin{equation}
  \label{eq-gamma}
\ga \leq 2 (\al + \be)
\end{equation}
and
\begin{equation}
  \label{eq-t}
\frac{1}{3} (\al + \be + \ga) \leq t.
\end{equation}
Condition (ii) reads in our case as
\[
\min \{\al + \be + \ga + t, \al + \be + t \} \geq t + \frac{2}{3} (\al +\be + \ga),
\]
which is equivalent to
\[
3t \geq 2 \ga - \al - \be.
\]
Using Inequality (\ref{eq-gamma}) one checks that the last condition is implied by Inequality (\ref{eq-t}). This completes the argument.
\end{proof}

Suppose we are given $1 \leq \al \leq \be$ and want to choose $\ga, t$ such that the syzygy bundle of $I$ is semistable. Then there is only a finite number of choices for $\ga$ since we must have $\be \leq \ga \leq 2 (\al + \be)$, whereas we have infinitely many choices for $t$ as the only condition is $t \geq \frac{1}{3} (\al + \be + \ga)$.

If the syzygy bundle of $I$ is not  semistable, then $R/I$ must have
the WLP (\cite{BK}, Theorem 3.3). Combining this with Corollary
\ref{mod 3} and Lemma \ref{lem-sst}, our computer experiments
suggests  the following characterization of the presence of the WLP
in characteristic zero.

\begin{conjecture} \label{level wlp conj}
Let $I \subset R = K[x,y,z]$ be a level Artinian almost complete
intersection, i.e., $I$ is of the form
\[
(x^{\alpha+t},
y^{\beta+t}, z^{\gamma+t}, x^\alpha y^\beta z^\gamma),
\]
where $t > 0$ and, without loss of generality, $0 \leq \alpha \leq
\beta \leq \gamma$. Assume that $K$ is algebraically closed field of
characteristic zero. Then:
\begin{itemize}
\item[(a)] $R/I$ has the WLP if any of the following conditions is
satisfied:
\begin{itemize}
  \item[(i)] $\alpha = 0$,
  \item[(ii)] $\al + \be + \ga$ is not divisible by 3,
  \item[(iii)] $\ga > 2 (\al + \be)$,
  \item[(iv)] $t < \frac{1}{3} (\al + \be + \ga)$.
\end{itemize}

\item[(b)] Assume that $1 \leq   \alpha \leq
\beta \leq \gamma \leq 2 (\al + \be)$,\; $\alpha + \beta +
\gamma \equiv 0 \ (\hbox{\rm mod } 3)$, and $t \geq \frac{1}{3}
(\al + \be + \ga)$. Then $R/I$ fails to have the WLP if and only
if $t$ is even and  either of the following two conditions is
satisfied:
\begin{itemize}
\item [(i)] $\al$ is even, $\al = \be$ and
$\ga - \al \equiv 3 \ (\hbox{\rm mod } 6)$,;

\item[(ii)] $\alpha$ is odd and

\hspace*{2cm} $\al = \be$ and $\ga - \al \equiv 0 \
 (\hbox{\rm mod } 6)$,

or

 \hspace*{2cm} $\be = \ga$ and $\ga - \al \equiv 0 \
 (\hbox{\rm mod } 3)$.

\end{itemize}
Furthermore, in all of the above cases,  the Hilbert function has
``twin peaks.''
\end{itemize}
\end{conjecture}

Note that part (a) is true by Lemma \ref{lem:variable-missing},
Corollary \ref{mod 3},  and Lemma \ref{lem-sst}. Part (b) will be
discussed in the following section.


\section{A proof of half of Conjecture \ref{level wlp conj}}
\label{half conj}

We are going to establish sufficiency of the numerical conditions
given in Conjecture \ref{level wlp conj}(b) for failure of the WLP.
First we establish the claim about the twin peaks of the Hilbert
function.

\begin{lemma} \label{twin peaks}
Consider the ideal $I_{\alpha, \beta, \gamma, t} =  (x^{\alpha+t},
y^{\beta+t}, z^{\gamma+t}, x^\alpha y^\beta z^\gamma)$ which (by
Proposition \ref{codim3}) defines a level algebra.
 Assume
 that  $ t \geq \max  \left \{
\frac{2 \gamma -\alpha - \beta}{3},  \frac{ \alpha + \beta +
\gamma}{3} \right  \}$ and that $0 < \alpha \leq \beta \leq \gamma$.
Assume furthermore that $\alpha + \beta + \gamma \equiv 0\
(\hbox{\rm mod } 3)$.  Then the values of the Hilbert function of
$R/I_{\alpha, \beta,\gamma,t}$ in degrees
$\frac{2(\alpha+\beta+\gamma)}{3} +t-2$ and
$\frac{2(\alpha+\beta+\gamma)}{3} +t-1$ are the same.
\end{lemma}

\begin{proof}
By Proposition \ref{codim3}, we know the minimal free resolution of  $R/I_{\alpha, \beta,\gamma,t}$, which we can use to compute the Hilbert function in any degree.  We first claim that in the specified degrees, this computation has no contribution from the last and the penultimate free modules in the resolution.  To do this, it is enough to check that the degree $\frac{2(\alpha+\beta+\gamma)}{3} +t-1$ component of any summand in the penultimate free module is zero.  The first three summands correspond to the observation that
\[
- \frac{\alpha}{3} - \frac{\beta}{3} - \frac{\gamma}{3} - 1 < 0.
\]
Since $\alpha \leq \beta \leq \gamma$ and $a = \alpha +t$, $b = \beta+t$, and $c = \gamma +t$, we have only to check that
\[
\frac{2(\alpha+\beta+\gamma)}{3} +t-1 - \alpha - t - \beta - t < 0.
\]
This is equivalent to the inequality on $t$ in the hypotheses.

Now rather than explicitly computing the Hilbert functions in the two degrees, it is enough to express them as linear combinations of binomial coefficients and show that the difference is zero, using the formula $\binom{p}{2} - \binom{p-1}{2} = p-1$. This is a routine computation.
\end{proof}

\begin{theorem} \label{WLP and M}
Consider the level algebra $R/I_{\alpha,\beta,\gamma,t} =
R/(x^{\alpha + t}, y^{\beta + t}, z^{\gamma + t}, x^\alpha y^\beta
z^\gamma )$.  We make the following assumptions:

\begin{itemize}

\item $\displaystyle 0 < \alpha \leq \beta \leq \gamma \le
2(\alpha + \beta)$\hbox{\rm ;}

\item $\displaystyle t \ge  \frac{\alpha + \beta +
\gamma}{3}$\hbox{\rm ;}

\item  $\displaystyle \alpha + \beta + \gamma \equiv 0 \   (\hbox{\rm mod } 3)$.

\end{itemize}

Then there is a square matrix, $M$, with integer entries, having the
following properties.

\begin{itemize}

\item[(a)] $M$ is a $\displaystyle \left ( t + \frac{\alpha + \beta - 2\gamma}{3} \right ) \times  \left ( t + \frac{\alpha + \beta - 2\gamma}{3} \right ) $ matrix.

\item[(b)] If $\det M \equiv 0 \ \hbox{\rm (mod } p)$, where $p$
is the characteristic of $K$, then $R/I_{\alpha,\beta,\gamma,t}$
fails to have the WLP.  This includes the possibility that $\det M
= 0$ as an integer.

\item[(c)] If $\det M \not \equiv  0  \ \hbox{\rm (mod } p)$ then
$R/I_{\alpha,\beta,\gamma,t} $ satisfies the WLP.

\end{itemize}

\end{theorem}

\begin{proof}
We note first that the second bullet in the hypotheses implies (using the first bullet) that the following inequality also holds:
\[
t >  \frac{2 \gamma -\alpha - \beta -3}{3}
\]
We will use this fact without comment in this proof.

Thanks to Lemma \ref{twin peaks}, the values of the Hilbert
function of $R/I_{\alpha,\beta,\gamma,t}$ in degrees
$\frac{2(\alpha+\beta+\gamma)}{3} +t-2$ and
$\frac{2(\alpha+\beta+\gamma)}{3} +t-1$ are the same.  Hence
thanks to Proposition \ref{gen wlp}, checking whether or not the
WLP holds is equivalent to checking whether multiplication by a
general linear form between these degrees is an isomorphism or
not.

We will use Proposition \ref{splitting and surj}.  Let $L$ be a general linear form, let $\bar R = R/(L) \cong K[x,y]$ and let $\bar I$ be the image of $I_{\alpha, \beta, \gamma, t}$ in $\bar R$.  Note that $\bar I \cong (x^{\alpha +t}, y^{\beta+t}, \ell^{\gamma +t}, x^\alpha y^\beta \ell^\gamma)$, where $\ell$ is the restriction to $\bar R$ of $z$, and  thanks to Lemma \ref{lem-L-element} we will take $\ell = x+y$.

Of course $x^{\alpha+t}, y^{\beta+t}$ is a regular sequence.  Hence it suffices to check whether or not  there is an element $F \in \bar R$ of degree
\[
f := \alpha + t + \beta + t - \left [ \frac{2(\alpha + \beta + \gamma)}{3} +t-1 +2 \right ] = \frac{\alpha + \beta-2 \gamma}{3} +t-1 = \frac{\alpha+\beta+\gamma}{3} - \gamma +t-1,
\]
non-zero modulo $(x^{\alpha +t}, y^{\beta +t})$, such that $F \cdot \bar I \subset (x^{\alpha +t}, y^{\beta +t})$.  The latter condition is equivalent to
\begin{equation} \label{two conditions}
F \cdot (x+y)^{\gamma +t} \in (x^{\alpha+t}, y^{\beta+t}) \hbox{ and } F \cdot x^\alpha y^\beta (x+y)^\gamma \in (x^{\alpha+t}, y^{\beta+t}).
\end{equation}

\bigskip

\noindent \underline{Claim}:
\begin{enumerate}

\item {\em $\gamma \geq 2(\alpha + \beta) $ if and only if $F \cdot (x+y)^{\gamma+t}$ is automatically in $(x^{\alpha+t} , y^{\beta+t})$.  }

\item  {\em  $\displaystyle t \leq \frac{\alpha + \beta + \gamma}{3}$ if and only if $F \cdot x^\alpha y^\beta (x+y)^\gamma$ is automatically in $(x^{\alpha+t} , y^{\beta+t})$.  }

\end{enumerate}

\noindent  Indeed, the first inequality is equivalent to $\deg ( F \cdot (x+y)^{\gamma +t} ) \geq (\alpha +t) + (\beta +t)-1$, so every term of $F \cdot (x+y)^{\gamma +t} $ is divisible by either $x^{\alpha+t}$ or $y^{\beta+t}$.  The second inequality is equivalent to $\deg ( F \cdot x^\alpha y^\beta (x+y)^\gamma ) \geq (\alpha +t) + (\beta +t)-1$.  This establishes the claim.  Thanks to our hypotheses, then, the conditions in (\ref{two conditions}) add constraints on the possibilities for $F$.  We want to count these constraints.

Let $F = \lambda_0 x^f + \lambda_1 x^{f-1}y + \lambda_2 x^{f-2}y^2 + \dots + \lambda_{f-2} x^2 y^{f-2} + \lambda_{f-1} x y^{f-1} + \lambda_f y^f.$  We now consider how many conditions (\ref{two conditions}) imposes on the $\lambda_i$.  Consider the first product, which has degree $\frac{\alpha +\beta+ \gamma}{3} + 2t -1$.  A typical term in $F \cdot (x+y)^{\gamma +t}$ is some scalar times $x^i y^j$, where $i+j = \frac{\alpha +\beta+ \gamma}{3} + 2t -1$.  The set of all pairs $(i,j)$ for which $x^iy^j$ is not in $(x^{\alpha +t}, y^{\beta+t})$ is
\begin{equation} \label{ij for first prod}
\{ (i,j) \} =  \left \{  \left (\alpha +t-1, \frac{\beta + \gamma - 2\alpha}{3} +t  \right ) , \dots, \left ( \frac{\alpha + \gamma -2\beta}{3} +t , \beta +t-1 \right )  \right \}
\end{equation}
Since each such term has to vanish, this imposes a total of $\frac{2\alpha + 2\beta -\gamma}{3}$ conditions on the $\lambda_i$.  Similarly, consider the second product, which has degree $\frac{4 \alpha + 4 \beta+\gamma}{3} +t-1$.  A typical term in $F \cdot x^\alpha y^\beta (x+y)^\gamma$ is some scalar times $x^i y^j$, where $i+j = \frac{4 \alpha + 4 \beta+\gamma}{3} +t-1$.  The set of all pairs $(i,j)$ for which $x^iy^j$ is not in the ideal $(x^{\alpha+t}, y^{\beta+t})$ is
\begin{equation} \label{ij for second prod}
\{ (i,j) \} = \left \{ \left ( \alpha + t -1, \frac{\alpha + 4\beta +\gamma}{3} \right ), \dots, \left ( \frac{4\alpha + \beta + \gamma}{3}, \beta+t-1 \right )         \right \}
\end{equation}
This imposes a total of $t - \frac{\alpha + \beta + \gamma}{3} $
conditions, since we need all of these terms to vanish. Combining,
we have a total of $t+ \frac{\alpha + \beta - 2\gamma}{3} = f+1$
conditions.  Since there are $f+1$ variables $\lambda_i$, the
coefficient matrix is the desired square matrix. Now it is clear
that $\det M = 0$ (regardless of the characteristic) if and only
of the corresponding homogeneous system has a non-trivial
solution, i.e. there is a polynomial $F$ as desired, if and only
if $R/I_{\alpha, \beta, \gamma, t}$ fails to have the WLP.
\end{proof}

We can specifically give the matrix described in the last result.

\begin{corollary} \label{specific M}
The matrix in Theorem \ref{WLP and M} has the form
\[
M =
\left [
\begin{array}{cccccccccccccccccccccc}
\binom{\gamma }{\frac{\alpha+\beta + \gamma }{3}} & \binom{\gamma
}{\frac{\alpha+\beta + \gamma }{3} -1} & \dots & \binom{\gamma
}{\gamma -t+2} & \binom{\gamma }{\gamma -t+1}
\\ \\
&& \vdots \\ \\
\binom{\gamma }{t-1} & \binom{\gamma }{t-2} & \dots &
\binom{\gamma }{\frac{2\gamma -\alpha -\beta}{3} +1}
 & \binom{\gamma }{\frac{2\gamma -\alpha -\beta}{3 }} \\ \\
\binom{\gamma+t}{t+\beta -1} & \binom{\gamma+t}{t+\beta -2} & \dots & \binom{\gamma+t}{\frac{2(\beta +\gamma)-\alpha}{3}+1} & \binom{\gamma+t}{\frac{2(\beta +\gamma)-\alpha}{3}} \\ \\
\binom{\gamma+t}{t+\beta -2} & \binom{\gamma+t}{t+\beta -3} &
\dots &
\binom{\gamma+t}{\frac{2(\beta +\gamma)-\alpha}{3}} & \binom{\gamma+t}{\frac{2(\beta +\gamma)-\alpha}{3}-1} \\ \\ && \vdots \\ \\
\binom{\gamma+t}{t+\frac{\beta +\gamma -2\alpha}{3}} &
\binom{\gamma+t}{t-1+\frac{\beta +\gamma -2\alpha}{3}} & \dots &
\binom{\gamma+t}{\gamma -\alpha+2} & \binom{\gamma+t}{\gamma
-\alpha+1}

\end{array}
\right ]
\]
\end{corollary}

\begin{proof}
This is a tedious computation, but is based entirely on the proof of Theorem \ref{WLP and M}.  The top ``half'' of the matrix corresponds to the second product in (\ref{two conditions}), and the bottom ``half'' of the matrix corresponds to the first product.  Each row in the top ``half''  corresponds to one ordered pair in (\ref{ij for second prod}), and each row in the bottom ``half'' corresponds to one ordered pair in (\ref{ij for first prod}).
\end{proof}

The following corollary establishes the sufficiency of the numerical conditions given in Conjecture \ref{level wlp conj}.

\begin{corollary} \label{half of conj}
Let $K$ be an arbitrary field and $R = K[x,y,z]$.   Consider the
ideal $I_{\alpha, \beta, \gamma, t} = (x^{\alpha + t}, y^{\beta +
t}, z^{\gamma + t}, x^\alpha y^\beta z^\gamma )$, where $1 \leq
\alpha \leq \beta \leq \gamma$. Assume that one of the following
three cases holds:

\begin{enumerate}
\item $(\alpha, \beta, \gamma,t) = (\alpha, \alpha, \alpha +3\lambda ,t)$ with $\alpha$ even, $\lambda $ odd, $t \ge \alpha +\lambda $ even and $1\le \lambda \le \alpha$;

\item $(\alpha, \beta, \gamma,t) = (\alpha, \alpha, \alpha + 6 \mu, t)$ with $\alpha$ odd, $t \geq \alpha + 2 \mu$ even, and  $0 \leq \mu \leq \frac{\alpha -1}{2}$; or

\item $(\alpha, \beta, \gamma,t) = (\alpha, \alpha + 3 \rho, \alpha + 3 \rho,t)$ with $\alpha$ odd, $t \geq \alpha + 2 \rho$ even, and $\rho \geq 0$.

\end{enumerate}
Then $R/I_{\alpha, \beta, \gamma, t}$  fails to have the WLP.
\end{corollary}

\begin{proof}
Possible after an extension of the base field, we may assume that
$K$ is an infinite field. One can verify quickly (using the
constraints on the invariants given in the theorem) that the
hypotheses of Theorem \ref{WLP and M} hold here in all three cases
(the parity is important in some instances).  Hence it is only a
matter of identifying $M$, via Corollary \ref{specific M}, and
checking that in all the cases mentioned, $\det M = 0$.   We first
consider Case~(1).

By applying Corollary \ref{specific M}, we obtain the $(t-2\lambda) \times (t-2\lambda)$ matrix $M$ below, corresponding to $t-2\lambda$ homogeneous equations in $t-2\lambda$ unknowns.
\[
M =
\left [
\begin{array}{cccccccccccccccccccccc}
\binom{\alpha+3\lambda}{\alpha+\lambda} & \binom{\alpha+3\lambda}{\alpha +\lambda -1} & \dots & \binom{\alpha +3\lambda}{\alpha-t+\lambda +2} & \binom{\alpha+3\lambda}{\alpha-t+\lambda +1} \\ \\
\binom{\alpha+3\lambda}{\alpha+\lambda +1} & \binom{\alpha+3\lambda}{\alpha+\lambda} & \dots & \binom{\alpha+3\lambda}{\alpha-t+\lambda +3} & \binom{\alpha+3\lambda}{\alpha-t+\lambda +2} \\
&& \vdots \\
\binom{\alpha+3\lambda}{t-1} & \binom{\alpha+3\lambda}{t-2} & \dots & \binom{\alpha+3\lambda}{2\lambda +1} & \binom{\alpha+3\lambda}{2\lambda} \\ \\
\binom{\alpha+t+3\lambda}{t+\lambda} & \binom{\alpha+t+3\lambda}{\lambda -1} & \dots & \binom{\alpha+t+3\lambda}{3\lambda+2} & \binom{\alpha+t+3\lambda}{3\lambda +1} \\ \\
\binom{\alpha+t+3\lambda}{t+\lambda +1} & \binom{\alpha+t+3\lambda}{t+\lambda} & \dots & \binom{\alpha+t+3\lambda}{3\lambda+3} & \binom{\alpha+t+3\lambda}{3\lambda+2} \\
&& \vdots \\
\binom{\alpha+t+3\lambda}{\alpha+t-1} & \binom{\alpha+t+3\lambda}{\alpha+t-2} & \dots & \binom{\alpha+t+3\lambda}{\alpha+2\lambda +1} & \binom{\alpha+t+3\lambda}{\alpha+2\lambda}

\end{array}
\right ]
\]
This system has a non-trivial solution (giving the existence of the desired form $F$) if and only if $M$ has determinant zero.

We will show that under our assumptions, this determinant is indeed zero.  Observe that if $M$ is flipped about the central vertical axis, and then the top portion and the bottom are (separately) flipped about their respective central horizontal axes, then we restore the matrix $M$.  Since $t-2\lambda$ is even, the first step can be accomplished with $\frac{t-2\lambda}{2}$ interchanges of columns.  The top portion contains $t-\lambda-\alpha$ rows and the bottom portion contains $\alpha-\lambda$ rows. Both are odd, so the second step can be done with $\frac{t-\lambda -1-\alpha}{2}$ interchanges and the last one with $\frac{\alpha-\lambda -1}{2}$ interchanges.  All together we have $t-2\lambda-1$ interchanges of rows/columns, which is an odd number.  Therefore $\det M = - \det M$, and so $\det M = 0$.

Case (2) is similar and is left to the reader.  Case (3), however, is somewhat different.  Now we will restrict to $K[y,z]$ rather than $K[x,y]$.  We obtain
\[
\bar I = ((y+z)^{\alpha+t}, y^{\alpha + 3 \rho +t}, z^{\alpha +3\rho +t}, (y+z)^{\alpha}\cdot y^{\alpha+3\rho} z^{\alpha +3\rho})
\]
and we have to check whether there is a form $F \in K[y,z]$ of degree $t + 2\rho-1$ (obtained after a short calculation) such that
\[
F \cdot (y+z)^{\alpha +t} \in (y^{\alpha +3\rho+t}, z^{\alpha + 3\rho +t}) \hbox{ and } F \cdot (y+z)^\alpha \cdot y^{\alpha +3\rho} z^{\alpha +3\rho} \in (y^{\alpha +3\rho+t}, z^{\alpha + 3\rho +t}).
\]
The calculations again follow the ideas of Theorem \ref{WLP and M}, and we obtain the following $(t + 2\rho) \times (t +2\rho)$ matrix of integers:

{\footnotesize
\[
\left [
\begin{array}{cccccccccccccccccccccccccccccc}
0 & 0 & \binom{\alpha }{0} & \binom{\alpha }{1} & \dots &  \dots &
\binom{\alpha }{\alpha} & 0 & 0 & \dots & 0 & \dots & 0
\\ \\
&& \vdots && \vdots & \vdots &&&& \vdots && \vdots\\ \\
0 & 0 & 0 & 0 & \binom{\alpha }{0} & \binom{\alpha }{1} & \dots &
\dots &  \dots & \binom{\alpha }{\alpha} & 0 & \dots & 0
  \\ \\

\binom{\alpha+t}{\alpha + \rho} & \binom{\alpha+t}{\alpha+\rho +1}
& \dots & \dots
 & \binom{\alpha+t}{\alpha +t} & 0 & \dots & \dots & \dots  & \dots & \dots & \dots & 0\\ \\
\binom{\alpha+t}{\alpha + \rho-1} & \binom{\alpha+t}{\alpha+\rho }
& \dots &  \dots &
\binom{\alpha+t}{\alpha +t-1} & \binom{\alpha+t}{\alpha +t} & 0 & \dots & \dots & \dots & \dots & \dots & 0\\ \\
&& \vdots &&& \vdots &&& \vdots && \vdots\\ \\
\binom{\alpha+t}{0} & \binom{\alpha+t}{1} & \dots & \dots
 & \dots  & \dots & \dots & \binom{\alpha+t}{\alpha+t} & 0 & \dots & \dots & \dots & 0\\ \\
 0 & \binom{\alpha+t}{0} & \dots & \dots
 & \dots  & \dots &  \dots & \binom{\alpha+t}{\alpha+t-1} & \binom{\alpha+t}{\alpha+t-1} & \dots &  \dots & \dots & 0\\ \\
&& \vdots &&& \vdots &&& \vdots &&& \vdots\\ \\
0 & 0 & 0 & \binom{\alpha+t}{0} & \dots & \dots & \dots & \dots
 & \dots & \dots & \dots & \dots & \binom{\alpha+t}{t-\rho} \\ \\

\end{array}
\right ]
\]
}
We remark that in the top portion (i.e.\ the submatrix where the binomial coefficients have $\alpha$ as the top entry), the first row has a sequence of $2\rho$ zeroes before the $\binom{\alpha}{0}$, and the last row (of the top portion) has a sequence of $2\rho$ zeroes after the $\binom{\alpha}{\alpha}$.  The top portion has $t-2\rho-\alpha$ rows, while the bottom portion (with binomial coefficients having $\alpha+t$ as top entry) has $\alpha+4\rho$ rows.
In the same way as before (using the fact that $t$ is even), it is easy to see that this matrix can be restored to itself with an odd number of row and column interchanges, and hence the determinant is zero.
\end{proof}

\begin{remark}
Notice that to check the surjectivity of the multiplication by a linear form from degree $d-1$ to degree $d$, we have to check whether or not $(R/(I,L))_d$ is zero.  Because of this, it is possible to obtain the result of Theorem \ref{half of conj} (including exactly the same matrix $M$) with a more direct computation, rather than using the liaison approach of Proposition \ref{splitting and surj}.  However, the computations seemed slightly more intricate, and we also felt that the existence of the form $F$ might have other interesting applications.
\end{remark}

\begin{corollary} \label{kkkaaa}
Consider the level algebras of the form $R/I$ with $I = (x^k, y^k,
z^k, x^\alpha y^\alpha z^\alpha)$, $\alpha$ odd and $k \geq 2
\alpha+1$ odd.   Then $R/I$ is level and fails to have the WLP.
\end{corollary}

\begin{example} \label{both cases}
Consider an ideal of the form $I = (x^{10}, y^{10}, z^{10}, x^3 y^3
z^3)$ (i.e. we relax the condition in Corollary \ref{kkkaaa} that
$k$ be odd).  Then $\det M = 78,\hbox{\hskip -.05cm}408 = 2^3 \cdot
3^4 \cdot 11^2$. One can check on a computer program (e.g. CoCoA
\cite{cocoa}) that in characteristic 2, 3 and 11, $R/I_{10}$ does
not have the WLP, while in characteristic 5, 7, 13, 17, \dots ,
$R/I$ does have the WLP  (as predicted by Theorem \ref{WLP and M}).
\end{example}

\begin{corollary}
 For any even socle degree there is a level monomial almost complete intersection which fails to have
 the WLP.
\end{corollary}

\begin{proof}
In Corollary \ref{kkkaaa}, simply consider the special cases $\alpha = 1$ and $\alpha = 3$.
>From Proposition \ref{codim3} we note that in the first case  the socle degree is $2k -2$ and in the second case the socle degree is $2k$.
\end{proof}

 The following is a natural question to ask  at this point:

\begin{question} \label{A(n) question}  Is there is a monomial level almost complete intersection (or indeed any almost complete intersection) in three variables with odd socle degree and
failing to have the WLP?
\end{question}

We now address this question in characteristic $p$.  We begin with a simple example:

\begin{example} \label{A(n) example}
Fix any prime $p$ and consider the complete  intersection
${\mathfrak a} = (x^p,y^p,z^p)$ in $R = K[x,y,z]$, where $K$ has
characteristic $p$.  Note that $R/{\mathfrak a}$ has socle degree
$3p-3$ and fails to have the WLP, since for a general linear form
$L$, $L^{p-1}$ is in the kernel of $(\times L)$.  (This was
observed for $p=2$ in \cite{HMNW}, Remark 2.9, and in \cite{MNZ3},
Remark 2.6, in arbitrary characteristic.)  Now consider the ideal
$I = ({\mathfrak a}, x^{p-1}y^{p-1}z^{p-1})$.  This clearly is an
almost complete intersection with socle degree $3p-4$ (an odd
number), is level, and still fails to have the WLP.

Notice that similar examples exist whenever the number of variables
is at least three.

\end{example}

This leads to the following refinement of our question:

\begin{question} For fixed characteristic $p$, what odd socle
degrees can occur for almost complete intersections without the
WLP?
\end{question}

Obviously quotients of the ideal $\mathfrak a$ cannot give us any
examples for socle degree $> 3p-3$, so we have to look to
different powers of the variables.  In \cite{HMNW} Remark 2.9 it
was observed that in characteristic 2 the ideal $(x^4,y^4,z^4)$
also fails to have the WLP.  We are led to consider other powers,
and we ask the following natural question.

\begin{question}
Given a prime $p$, consider  ideals $(x^k, y^k, z^k)$ in
characteristic $p$.  For which values of $k$ does
$R/(x^k,y^k,z^k)$ fail to have the WLP?
\end{question}

In characteristic zero, on the other hand, the situation  seems to
be different.  We conjecture that the converse of Corollary
\ref{half of conj} holds, which is the only missing piece in
characterizing the level monomial ideals in $K[x,y,z]$ that fail to
have the WLP and establishing Conjecture \ref{level wlp conj}.

\begin{conjecture} \label{conj with bad wording}
Assume that $K$ has characteristic zero.  Using the notation of
Corollary \ref{half of conj}, if $R/I_{\alpha, \beta, \gamma, t}$
fails to have the WLP then one of cases (1), (2) or (3) holds.
\end{conjecture}

\begin{remark}
  \label{rem-what-remains}
According to Lemma \ref{twin peaks}, the Hilbert function of
$R/I_{\alpha, \beta, \gamma, t}$ agrees in degrees $s =
\frac{2(\alpha+\beta+\gamma)}{3} +t-2$ and $s+1$. Hence, by
Proposition \ref{gen wlp},  Conjecture \ref{conj with bad wording}
and thus also Conjecture \ref{level wlp conj} is proven, if one
shows that the multiplication map $(R/I_{\alpha, \beta, \gamma,
t})_s \to (R/I_{\alpha, \beta, \gamma, t})_{s+1}$ by $x+y+z$ is
injective or surjective, provided the conditions in (1), (2) or (3)
of Corollary \ref{half of conj} all fail to hold.
\end{remark}

\begin{remark}
We have been focusing on monomial complete intersections.   A
natural question is whether a ``general'' complete intersection of
height three in characteristic $p$ necessarily has the WLP even
when the monomial one does not.  If the field $K$ is at least
infinite, this is answered in the affirmative by the main result
of \cite{anick}.
\end{remark}


\section{Final Comments} \label{final comments}

In Section \ref{arb codim} we saw that for the ideal $I = (
x_1^r,\dots,x_r^r,x_1\cdots x_r ) \subset R = K[x_1,\dots,x_r]$,
the corresponding algebra $R/I$ fails to  have the WLP.  On the
other hand, we saw in Section \ref{almost monomial} that making a
very slight change to even one of the generators gave an algebra
with the same Hilbert function, but possessing the WLP.  In this
section we analyze related phenomena and pose related questions.

\begin{example} \label{codim 4}
Stanley \cite{stanley} and Watanabe \cite{watanabe}  showed that a
monomial complete intersection satisfies the {\em Strong Lefschetz
Property (SLP)}.  This property is a generalization of the WLP,
and says that if $L$ is a general linear form then for any $i$ and
any $d$, the multiplication $\times L^d : (R/I)_i \rightarrow
(R/I)_{i+d}$ has maximal rank.  This means that the algebra
$R/(I,L^d)$ has the ``expected'' Hilbert function.  By
semicontinuity, the same is true when $L^d$ is replaced by a
general form of degree $d$.

It is of interest to find behavior of $R/I$ that distinguishes multiplication by $L^d$ from that by a general form of degree $d$.  We have found such a phenomenon experimentally on CoCoA \cite{cocoa}, although we have not given a theoretical justification.  Let $R = K[x_1,x_2,x_3,x_4]$ and let $L \in R$ be a general linear form.  Consider the ideals
\[
\begin{array}{rcl}
I_N & = & (x_1^N, x_2^N, x_3^N, x_4^N, L^N); \\ \\
J_N & = & (x_1^N, x_2^N, x_3^N, x_4^N, G)
\end{array}
\]
where $G$ is a general form of degree $N$.
  By the above-cited result, $R/J_1$ and $R/J_2$ have the same Hilbert function, and it can be checked that in fact they have the same minimal free resolution.  However, these algebras often have different behavior with respect to the
   WLP and with respect to minimal free resolutions!
More precisely, we have the following experimental data, which we computed in \cocoa \ over the rational numbers.

\bigskip

\begin{center}
\begin{tabular}{c|c|c|ccccccccccccccccccccccccc}
$N$ & $I_N$ has the WLP? & $J_N$ has the WLP? & Same resolution? \\
\hline
2 & yes & yes & yes\\
3 & no & yes & no \\
4 & no & yes & yes \\
5 & no & yes & yes \\
6 & no & yes & yes \\
7 & no & yes & no \\
8 & no & yes & yes \\
9 & no & yes & yes \\
10 & no & yes & yes \\
11 & no & yes & no \\
12 & no & yes & ?
\end{tabular}
\end{center}

\bigskip

\end{example}

\begin{question} \label{list of questions}

We end by posing some natural questions that remain to be  addressed
(in addition to the conjectures posed earlier).

\begin{enumerate}


\item Have all the classes of algebras studied in this paper a
unimodal  Hilbert function?



\item For $r = 3$ we have the following questions.
 (Here ``ACI'' means ``almost complete intersection.")

\begin{enumerate}
\item Are monomial ideals the only ACI's that fail to have the
WLP?

\item Are there level ACI's without the WLP and having odd socle
degree?

\item Do there exist ACI's that fail to have  the WLP and are not
level?

\item Has every  monomial ideal (not necessarily ACI) that is
level of type two the  WLP?
\end{enumerate}



\item Have all ACI's a unimodal Hilbert function?


\end{enumerate}
\end{question}

\end{document}